\documentclass{amsart}
\usepackage{amssymb,euscript,amsmath, mathrsfs}
\usepackage[dvips]{graphicx}
\usepackage[dvips]{color}
\usepackage{epsfig,pstricks}

\newcounter{ENUM}
\newcommand{\itm}{\item}
\newenvironment{ilist}{\renewcommand{\theENUM}{\roman{ENUM}}\renewcommand{\itm}{\addtocounter{ENUM}{1}\item[(\theENUM)]}\begin{itemize}\setcounter{ENUM}{0}}{\end{itemize}}
\newenvironment{alist}[1][0]{\renewcommand{\theENUM}{\alph{ENUM}}\renewcommand{\itm}{\addtocounter{ENUM}{1}\item[\theENUM)]}\begin{itemize}\setcounter{ENUM}{#1}}{\end{itemize}}

\newtheorem{thm}{Theorem}[section]
\newtheorem{prop}[thm]{Proposition}
\newtheorem{lem}[thm]{Lemma}
\newtheorem{cor}[thm]{Corollary}

\theoremstyle{definition}
\newtheorem{defn}[thm]{Definition}

\newtheorem{ex}[thm]{Example}

\theoremstyle{remark}

\newtheorem{rem}[thm]{Remark}

\numberwithin{equation}{section}

\def\N{{\mathbb N}}
\def\Z{{\mathbb Z}}

\def\C{{\bf C}}

\def\supp{\operatorname{supp}}
\def\fac{\operatorname{fac}}

\def\Fac{\mathcal Fac}

\def\cG{\mathcal G}
\def\L{\mathcal L}
\def\l{\mathfrak l}
\def\MR{\mathcal{MR}}
\def\LMR{\mathcal{LMR}}
\def\R{\mathcal R}

\subjclass[2010]{05A15}

\begin{document}
\title{Factorizations of cycles and multi-noded rooted trees}
\author{Rosena R. X. Du}
\address{Rosena R. X. Du, Department of Mathematics, Shanghai Key Laboratory of Pure Mathematics and Mathematical Practice, East China Normal University, 500 Dongchuan Rd, Shanghai, 200041 P. R. China. Tel: 86 21 5434 2646-431, Fax: 86 21 5434 2633.}
\email{rxdu@math.ecnu.edu.cn}

\thanks{Rosena R. X. Du is partially supported by the National Science Foundation of China under Grant No. 10801053, Shanghai Rising-Star Program (No. 10QA1401900), and the Fundamental Research Funds for the Central Universities.}

\author{Fu Liu}
\thanks{Fu Liu is partially supported by the National Security Agency under Grant No. H98230-09-1-0029, and the National Science Foundation of China under Grant No. 10801053.}
\address{Fu Liu, Department of Mathematics,
University of California, Davis,
One Shields Avenue,
Davis, CA 95616 USA.}
\email{fuliu@math.ucdavis.edu}

\begin{abstract}
  In this paper, we study factorizations of cycles. The main result is that 
  under certain condition, the number of ways to factor a $d$-cycle into a product of cycles of prescribed lengths is $d^{r-2}.$
 To prove our result, we first define a new class of combinatorial objects, multi-noded rooted trees, which generalize rooted trees. We find the cardinality of this new class which with proper parameters is exactly $d^{r-2}.$ The main part of this paper is the proof that there is a bijection from factorizations of a $d$-cycle to multi-noded rooted trees via factorization graphs. This implies the desired formula.

  The factorization problem we consider has its origin in geometry, and is related to the study of a special family of Hurwitz numbers: pure-cycle Hurwitz numbers.
  Via the standard translation of Hurwitz numbers into group theory, our main result is equivalent to the following: when the genus is $0$ and one of the ramification indices is $d,$ the degree of the covers, the pure-cycle Hurwitz number is $d^{r-3},$ where $r$ is the number of branch points.

\end{abstract}

\maketitle

\section{Introduction}
Suppose $d$ is a positive integer and $\lambda = (\lambda_1, \lambda_2, \dots, \lambda_\ell)$ is a partition of $d.$ We say a permutation $\sigma \in S_d$ has {\it cycle type} $\lambda$ if $\lambda_1, \dots, \lambda_\ell$ are the lengths of the cycles in the cycle decomposition of $\sigma.$ We call a permutation $\sigma \in S_d$ an $e$-cycle if its cycle type is $(e,1, \dots, 1)$ for some $e \geq 2.$ Given a permutation $\sigma$ of cycle type $\lambda,$ we define its {\it index} as $\iota(\sigma) = \iota(\lambda) = \sum_{i} (\lambda_i-1).$ For any $e$-cycle $\tau,$ we define the {\it support} of $\tau$, denoted by $\supp(\tau)$, to be the set of $e$ elements that appear in $\tau.$

\begin{defn}\label{defn:fac}
Fix a $d$-cycle $\tau.$  We say $(\sigma_1, \dots, \sigma_{r-1})$ is a {\it factorization of $\tau$} if the following conditions are satisfied:
\begin{alist}
\itm For each $i,$ $\sigma_i$ is a cycle in $S_{\supp(\tau)};$
\itm $\sigma_1 \cdots \sigma_{r-1} = \tau.$
\end{alist}

If further for each $i,$ $\sigma_i$ is an $e_i$-cycle, we say $(\sigma_1, \dots, \sigma_{r-1})$ is a {\it factorization of $\tau$ of type $(e_1, \dots, e_{r-1})$}.

We denote by $\Fac(d,r,\tau; e_1, \dots, e_{r-1})$ the set of all the factorizations of $\tau$ of type $(e_1, \dots, e_{r-1})$ and $\fac(d,r,\tau; e_1, \dots, e_{r-1})$ the cardinality of $\Fac(d,r,\tau;$ $e_1,$ $\dots,$ $e_{r-1})$.
\end{defn}

Clearly, the number of factorizations is independent of the choice of $\tau$, so we often omit $\tau$ and just write $\fac(d,r;e_1,\dots,e_{r-1})$.

Our main result is the following theorem.

\begin{thm}\label{main2}
Suppose $\sum_{i=1}^{r-1} (e_i - 1) = d-1.$ Then
\begin{equation}
\fac(d,r; e_1, \dots, e_{r-1}) = d^{r-2}.
\end{equation}
\end{thm}

We remark that if $e_1 = \cdots = e_{r-1} = 2$, then $d = r$ and $\fac(d,r; e_1, \dots, e_{r-1})$ counts the number of factorizations of a $d$-cycle into $d-1$ transpositions. According to our theorem, this number is
\begin{equation}\label{factran}
\fac(d,d; 2, \dots, 2) = d^{d-2}.
\end{equation}
Note that $d^{d-2}$ also counts the number of labeled trees with $d$ vertices. Different bijective proofs of \eqref{factran} were given by D\'{e}nes \cite{Denes}, Moszkowski \cite{Moszkowski}, Goulden-Pepper \cite{GouldenPepper} and Goulden-Yong \cite{Treelike}.

\subsection*{Geometric background}
The factorization problem we consider arises from geometry. In this part, we briefly discuss the connection between the enumeration of factorizations of a cycle to counting a special case of Hurwtiz numbers, and conclude with a result on ``pure-cycle Hurwtiz numbers'' that is equivalent to Theorem \ref{main2}. The contents discussed here are irrelevant to the rest of the paper. The reader should feel free to skip it.

Hurwitz numbers count the number of connected branched covers of the projective line with specified ramification.
More precisely, the {\it Hurwitz number $h(d,r,g; \lambda^1, \dots, \lambda^r)$} counts the number of connected genus-$g$ covers of the projective line of degree $d$ with $r$ branch points where the monodromy over the $i$th branch point has cycle type $\lambda^i$. If a cover has non-trivial automorphisms, we divide by the size of its automorphism group.
According to the Riemann-Hurwitz formula, a branched cover satisfies
\begin{equation}\label{RHF}
\sum_{i=1}^r \iota(\lambda^i) = 2d-2+2g.
\end{equation}
Therefore, we are only interested in data $(d,r,g; \lambda^1, \dots, \lambda^r)$ that satisfies the above formula.

There is a group-theoretic description of Hurwitz numbers.

\begin{defn}
  Suppose $d \ge 1,$ $g \ge 0,$ $r \ge 0,$ and $\lambda^1, \dots, \lambda^r$ are partitions of $d$ satisfying \eqref{RHF}. A {\it Hurwitz factorization of type $(d,r,g; \lambda^1, \dots, \lambda^r)$}
is a tuple $(\sigma_1, \dots, \sigma_r)$ satisfying:
\begin{alist}
\itm $\sigma_i \in S_d$ has cycle type $\lambda^i;$
\itm $\sigma_1 \cdots \sigma_r = 1;$
\itm the $\sigma_i$'s generate a transitive subgroup of $S_d.$
\end{alist}
The {\it Hurwitz number $h(d,r,g; \lambda^1, \dots, \lambda^r)$} is the number of Hurwitz factorizations divided by $d!$.
\end{defn}

There has been a lot of work on Hurwitz numbers. Most of it has studied situations where all but one or two branch points are simple; i.e., all but one or two $\lambda^i$'s have the form $(2,1,\dots, 1).$ Hurwitz \cite{hurwitz} and Goulden-Jackson \cite{GouldenJakson} showed that if $\lambda^1 = \cdots = \lambda^{r-1} = (2,1,\dots,1)$ and $\lambda^r = (\tau_1, \dots, \tau_n),$ then
\begin{equation}\label{simple}
h(d,r,0; \lambda^1, \dots, \lambda^r) = \frac{(r-1)! d^{n-3} \prod_{i=1}^n \tau_i^{\tau_i}/\tau_i!}{m_1! m_2! \cdots m_d!},
\end{equation}
where $m_i$ is the number of $i$'s in $\lambda^r$ for $1 \le i \le d.$

Another special case of the Hurwitz numbers that has been studied is the the {\it pure-cycle} Hurwitz numbers. We say a Hurwitz number is {\it pure-cycle} if each $\lambda^i$ is of the form $(e_i, 1, \dots, 1)$ for some integer $e_i \ge 2.$ In other words, a pure-cycle Hurwitz number counts the number of genus-$g$ covers of the projective line of degree $d$ with $r$ branch points where there is only one ramification point over each branch point, with ramification index $e_i.$ In this situation, we will abbreviate our notation for the Hurwitz number to $h(d,r,g; e_1, \dots, e_r).$ Pure-cycle Hurwitz numbers were first studied in \cite{LiuOsserman}. The authors showed that
\begin{equation}\label{4point}
h(d,4,0; e_1, e_2, e_3, e_4) = \min\{e_i(d+1-e_i)\}.
\end{equation}

The number of factorizations of a $d$-cycle we consider is in fact related to a special case of pure-cycle Hurwitz numbers when one of the $e_i$'s is $d.$ Since the order of $e_i$'s does not change the Hurwitz number, without loss of generality, we can assume $e_r =d.$  Note that by the Riemann-Hurwitz formula \eqref{RHF}, we must have $\sum_{i=1}^{r} (e_i-1) = \sum_{i=1}^{r-1} (e_i -1) + (d-1) = 2d - 2 + 2 \cdot g.$ Hence, we require
\begin{equation}\label{RHF2}
\sum_{i=1}^{r-1} (e_i - 1) = d-1+2g.
\end{equation}
We remark that since $\sigma_r$ is a $d$-cycle, condition c) in the definition of Hurwitz factorization is automatically satisfied. Thus, to verify whether one tuple $(\sigma_1, \dots, \sigma_r)$ is a Hurwitz factorization of type $(d,r,g; e_1, \dots, e_{r-1}, d),$ we only need to check whether a) and b) are satisfied, which are precisely corresponding to conditions a) and b) in Defintion \ref{defn:fac}. Furthermore, in $S_d,$ there are $(d-1)!$ permutations that are $d$-cycles. Thus, the number of Hurwitz factorizations of type $(d,r,g; e_1, \dots, e_{r-1}, d)$ is $(d-1)! \fac(d,r;e_1\dots, e_{r-1}).$ Hence, assuming \eqref{RHF2}, we have that $$h(d,r,g; e_1, \dots, e_{r-1}, d) = \frac{1}{d} \fac(d,r; e_1, \dots, e_{r-1}).$$
Finally, we focus on the cases where $g=0.$ Then condition \eqref{RHF2} becomes
\begin{equation}\label{req}
\sum_{i=1}^{r-1} (e_i - 1) = d-1.
\end{equation}
Therefore, we conclude that Theorem \ref{main2} is equivalent to the following theorem.

\begin{thm}\label{main1}
Suppose $\sum_{i=1}^{r-1} (e_i - 1) = d-1.$ Then
\begin{equation}
h(d,r,0; e_1, \dots, e_{r-1}, d) = d^{r-3}.
\end{equation}
\end{thm}

One checks that Theorem \ref{main1} agrees with \eqref{simple} and \eqref{4point} in the corresponding special cases.

\subsection*{Other related work and organization of the paper}
Different but equivalent versions of Theorem \ref{main2} have previously been studied. Given nonnegative integers $n_2,\dots,n_d$, we say a factorization $(\sigma_1,\dots,\sigma_{r-1})$ of a $d$-cycle is {\it of cycle index $(n_2,n_3,\dots,n_d)$} if there are $n_m$ $m$-cycles among $\sigma_1,\dots, \sigma_{r-1}$ for any $2 \le m \le d$. Note that with this definition, the condition \eqref{req} translates to
\begin{equation}\label{req2}
\sum_{m=1}^d (m-1) n_m = d-1.
\end{equation}
Springer \cite{Springer} and Irving \cite{Irving} showed that assuming \eqref{req2}, the number of factorizations $(\sigma_1,\dots,\sigma_{r-1})$ of a $d$-cycle of cycle index $(n_2,\dots,n_d)$ is given by
\begin{equation}\label{sym}
\displaystyle d^{r-2} \frac{(r-1)!}{\prod_{m=2}^d n_m!}.
\end{equation}
Since the factorization number $\fac(d,r;e_1,\dots,e_{r-1})$ we consider is invariant under order of $e_i$'s, we see that Theorem \ref{main2} is equivalent to their result.
Springer \cite{Springer} proved the result by symmetrizing the problem further. He gave a bijection between factorizations of cycle index $(n_2,\dots,n_d)$ of {\it all} $d$-cycles in $S_d$ and doubly-labeled oriented cacti preserving cycle lengths, then showed the latter class of combinatorial objects has cardinality $(d-1)!$ times \eqref{sym}.
Irving's proof \cite{Irving} is based on a bijection between factorizations of cycle index $(n_2,\dots,n_d)$ of a {\it fixed} $d$-cycle and proper polymaps. (Irving's polymap is a generalization of the oriented cactus in \cite{Springer}. It can be used in general factorization problems without the restriction that each $\sigma_i$ has to be a cycle.)

Goulden-Jackson \cite{GouldenJacksonEJC} give a more general definition of factorizations of a $d$-cycle where they allow $\sigma_i$ to be any cycle type, that is, $\sigma_i$ does not have to be a cycle. They proved a more general result than Theorem \ref{main2}: Suppose $\lambda^1, \dots, \lambda^{m}$ are partitions of $d,$ where $\lambda^i$ consists of $\lambda_{j}^i$ $j$'s. Then the number of factorizations $(\sigma_1, \dots, \sigma_{m})$ of a fixed $d$-cycle of type $(\lambda^1, \dots, \lambda^{m})$ is given by
\begin{equation}\label{equ:general}
	d^{m-1} \frac{\prod_{i=1}^{m} (\sum_j \lambda_{j}^i -1)!}{\prod_{i=1}^{m} \prod_j \lambda_{j}^i}.
\end{equation}
Goulden-Jackson gave a bijection between the factorizations of a $d$-cycle of a given type to the plane-edge rooted $m$-cacti on $d$ $m$-gons with corresponding vertex distribution, and use a generating function argument to prove that the number of cacti is given by \eqref{equ:general}. Lando-Zvonkine gave a completely different proof for Goulden-Jackson's result in \cite{LandoZvonkine}.

The proofs given in \cite{Springer, Irving} can be considered as symmetrized bijective proofs of Theorem \ref{main2}. For the proofs for the more general result \eqref{equ:general}, techniques of Goulden-Jackson \cite{GouldenJacksonEJC} involve calculations with generating functions, and the paper of Lando-Zvonkine \cite{LandoZvonkine} uses geometric arguments.
In contrast, our techniques give a direct ``de-symmetrized'' bijective proof for Theorem \ref{main2}.
In order to do that, we first construct a new class of combinatorial objects called {\it multi-noded rooted trees}, show that (with proper parameters) it has cardinality $d^{r-2},$ and then give a bijection between factorizations of a $d$-cycle and multi-noded rooted trees.

The plan of this article is as follows: In Section \ref{sec:multi}, we define multi-noded rooted trees and find its cardinality. In Section \ref{sec:bigraph}, we associate to each factorization a bipartite graph, which we call factorization graph, and show that this association is injective. In Section \ref{sec:g2t}, we define a map from factorization graphs to multi-noded rooted trees, and state in Theorem \ref{thm:bijg2t} that it is a bijection. Assuming the theorem, we conclude Theorems \ref{main2} and \ref{main1}. In Section \ref{sec:charGraph}, we give characterizations of factorization graphs. Using this characterization, we complete the proof of Theorem \ref{thm:bijg2t} in Section \ref{sec:proof}.

Finally, we remark that a few results presented in Sections \ref{sec:bigraph} and \ref{sec:charGraph} have analogous or equivalent forms in the literature. However, to make our papers self-contained and accessible to readers without previous knowledge of geometric background of Hurwitz numbers, we include our proofs, which are purely combinatorial and only based on the group-theoretic description of Hurwitz numbers.

\subsection*{Acknowledgements}
We would like to thank Brian Osserman for providing data on pure-cycle Hurwitz numbers to us. We are also grateful to Richard Stanley who pointed out to us the reference \cite{Springer} and sent us a copy of it.

\section{Multi-noded rooted trees}\label{sec:multi}

We assume the readers are familiar with basic terminology in graph theory as presented in the appendix of \cite{stanleyec1}. We will review briefly the terms that will be used in this paper.

Recall that a graph is a pair $(V, E)$ where $V$ is the vertex set and $E \subseteq {V \choose 2}$ is the edge set of the graph.
A {\it tree} is an acyclic graph, and a {\it rooted tree} is a tree with a special vertex, which we call the {\it root} of the given tree.
Given a rooted tree $T,$ let $e = \{v,w\}$ be an edge of $T.$ If $v$ is closer to the root than $w,$ we call $v$ the {\it parent} of $w$ and $w$ a {\it child} of $v;$
we also call $v$ the {\it parent end} of $e$ and $w$ the {\it child end} of $e.$

We usually draw a rooted tree with its root at the top, put each child below the parent, and represent the vertices of the tree by distinct integers, i.e.,$V \subseteq \Z$.
In this paper, we always represent roots with the number $0.$  See Figure \ref{fig_prufer} for examples of rooted trees.

Suppose $S \subseteq \Z$ is a set of $n$ elements and $0 \not\in S.$ Let $\R_S$ be the set of rooted trees with vertex set $S \cup \{0\}$ and rooted at $0$. It is well-known that
\begin{equation}\label{RTcard}
|\R_S| = (n+1)^{n-1}.
\end{equation}

In this section, we will introduce a new class of combinatorial objects, called {\em multi-noded rooted trees}, which generalize $\R_S,$ and we will find its cardinality, which is exactly $d^{r-2}$ if we choose the right parameters.

Throughout this section, we assume $S=\{s_1 < s_2 < \cdots < s_{n}\}$ is a set of $n$ integers disjoint from $\{0\}.$

\begin{defn}
Suppose $f_0, f_1, \dots, f_{n}$ are positive integers. We say $M =(T, \beta)$ is a {\it multi-noded rooted tree} on $S \cup \{0\}$ of vertex data $(f_0, f_1, \dots, f_{n})$ if $T  = (S \cup \{0\}, E)$ is a rooted tree in $\R_S$ and $\beta: E \to \N$ is a function satisfying that for any edge $e \in E,$ if $s_i$ is the parent end of $e$, then $\beta(e) \in \{1, 2,\dots, f_i\}$.

We define $\MR_S(f_0, f_1, \dots, f_{n})$ to be the set of all multi-noded rooted trees on $S \cup \{0\}$ of vertex data $(f_0, f_1, \dots, f_{n}).$

We call the simple graph with one vertex and no edges {\it the trivial tree}.
\end{defn}

\subsection*{Graph representations of multi-noded rooted trees}
We give two ways to represent a multi-noded rooted tree $M = (T, \beta)$ graphically. The first way is to draw the rooted tree $T$ and then label each edge $e$ with $\beta(e).$ We call this the {\it edge-labeled representation} of $M.$

The second method is to draw a graph with {\it multi-noded} vertices: Given any positive integer $f,$ an {\it $f$-noded vertex} is a picture of $f$ nodes in a horizontal line and grouped together by a circle. (Note that the nodes in an $f$-noded vertex are considered to be ordered.) A {\it multi-noded} vertex is an $f$-noded vertex for some $f \in \N.$  With this definition, we can draw $M = (T, \beta)$ in the following way:
\begin{enumerate}
\item For each $0 \le i \le n,$ we draw an $f_i$-noded vertex which is labeled by $s_i.$ These $n+1$ multi-noded vertices are the vertices of $M.$
\item For any edge $e = \{s_i,s_j\}$ of $T$ with $s_i$ being the parent end of $e,$ we connect the multi-noded vertex $s_j$ to the $\beta(e)$-th node in vertex $s_i.$ These are the edges of $M.$
\end{enumerate}
We call this the {\it multi-noded representation} of $M.$

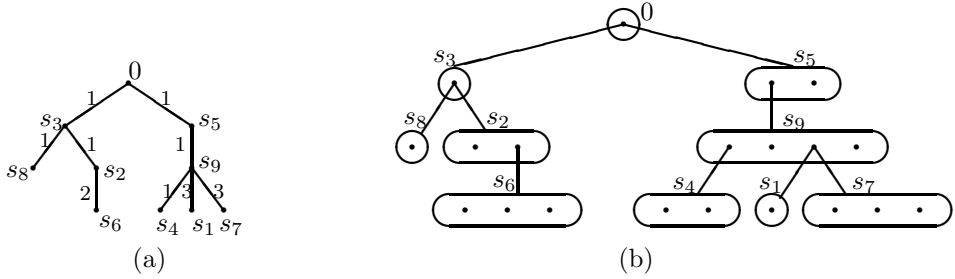
\begin{figure}[h,t]
\begin{center}
\setlength{\unitlength}{8pt} \thicklines

\begin{picture}(44,12)
\put(0,0){
\put(6,0){\makebox(0,0)[b]{(a)}}
\begin{picture}(14,8)

\put(5,9){\circle*{0.3}}\put(5,9.2){$0$}
\put(5,9){\line(-3,-2){3}}
\put(5,9){\line(3,-2){3}}

\put(2,7){\circle*{0.3}}\put(0.8,7){$s_3$}
\put(2,7){\line(-3,-4){1.5}}
\put(2,7){\line(3,-4){1.5}}

\put(0.5,5){\circle*{0.3}}\put(-0.8,4.6){$s_8$}
\put(3.5,5){\circle*{0.3}}\put(3.8,4.6){$s_2$}
\put(3.5,5){\line(0,-1){2}}
\put(3.5,3){\circle*{0.3}}\put(3.6,2.2){$s_6$}

\multiput(8,7)(0,-2){3}{\circle*{0.3}}
\put(8.3,7){$s_5$}
\put(8.3,5){$s_9$}
\put(8,2){$s_1$}\put(6.3,2){$s_4$}\put(9.3,2){$s_7$}
\put(8,7){\line(0,-1){4}}
\put(8,5){\line(-3,-4){1.5}}
\put(6.5,3){\circle*{0.3}}
\put(8,5){\line(3,-4){1.5}}
\put(9.5,3){\circle*{0.3}}
\footnotesize
\put(3,8){1}
\put(6.5,8){1}
\put(0.8,6){1}
\put(3,5.8){1}
\put(2.7,3.5){2}
\put(7.2,5.8){1}
\put(6.6,3.6){1}
\put(7.5,3.6){3}
\put(9,3.6){3}
\end{picture}}

\put(29,0){\makebox(0,0)[b]{(b)}}
\put(18,2){
\begin{picture}(26,8)
\put(10,9.8){\circle*{0.3}}\put(10.8,10){$0$}
\put(10,9.8){\oval(1.5,1.5)}
\put(10,9.8){\line(-4,-1){8}}
\put(10,9.8){\line(4,-1){8}}

\put(2,7){\circle*{0.3}}\put(1,8){$s_3$}
\put(2,7){\line(-2,-3){1.55}}
\put(2,7){\line(2,-3){1.5}}
\put(2,7){\oval(1.5,1.5)}

\put(0,4){\circle*{0.3}}\put(-0.4,5.1){$s_8$}
\put(0,4){\oval(1.5,1.5)}

\multiput(3,4)(2,0){2}{\circle*{0.3}}
\put(3.5,5.1){$s_2$}
\put(4,4){\oval(5,1.5)}
\put(5,4){\line(0,-1){2.25}}

\multiput(2.5,1)(2,0){3}{\circle*{0.3}}
\put(3.8,2.1){$s_6$}
\put(4.5,1){\oval(7,1.5)}

\multiput(17,7)(2,0){2}{\circle*{0.3}}\put(18,8){$s_5$}
\put(18,7){\oval(4.5,1.5)}
\put(17,7){\line(0,-1){2.25}}
\put(15,4){\line(-2,-3){1.5}}
\multiput(15,4)(2,0){4}{\circle*{0.3}}\put(18,4){\oval(9,1.5)}
\put(17.5,5){$s_9$}

\multiput(12,1)(2,0){2}{\circle*{0.3}}\put(13,1){\oval(5,1.5)}
\put(12.3,2){$s_4$}

\put(19,4){\line(-2,-3){1.6}}
\put(19,4){\line(2,-3){1.5}}

\put(17,1){\circle*{0.3}}
\put(17,1){\oval(1.5,1.5)}
\put(16.4,2){$s_1$}

\multiput(20,1)(2,0){3}{\circle*{0.3}} \put(22,1){\oval(7,1.5)}
\put(20.8,2){$s_7$}

\end{picture}}
\end{picture}
\caption{Two representations of a multi-noded rooted tree.}\label{tworep}
\end{center}
\end{figure}

\begin{ex}\label{ex_multi-noded}
Let $T=T_1$ as shown in Figure \ref{fig_prufer}. Suppose $M=(T,\beta)$ is the multi-noded rooted tree of vertex data $(1,1,2,1,2,2,3,3,1,4)$ and $\beta(\{0,s_3\})=1$, $\beta(\{0,s_5\})=1$, $\beta(\{s_3,s_8\})=1$, $\beta(\{s_3,s_2\})=1$, $\beta(\{s_5,s_9\})=1$, $\beta(\{s_2,s_6\})=2$, $\beta(\{s_9,s_4\})=1$, $\beta(\{s_9,s_1\})=3$, $\beta(\{s_9,s_7\})=3$. The two representations of $M$ are shown in Figure \ref{tworep}. Graph (a) is the edge-labeled representation and graph (b) is the multi-noded representation.
\end{ex}

\begin{rem}
We remark that each of the two representations has its own advantage. The edge-labeled representation does not involve new combinatorial structure. We will use it to find the cardinality of $\MR_S(f_0, f_1, \dots, f_n).$ The multi-noded representation contains the information of the vertex data while the edge-labeled representation does not. For example, graph (a) in Figure \ref{tworep} could be the graph of a multi-noded rooted tree of vertex data $(1,1,2,1,1,1,1,1,1,3)$ or anything bigger, but graph (b) in Figure \ref{tworep} can only be associated with vertex data $(1,1,2,1,2,2,3,3,1,4).$ The multi-noded representation will be used in a bijection we construct in Section \ref{sec:g2t}.
\end{rem}

\begin{prop}\label{MRcard}
The cardinality of $\MR_S(f_0, f_1, \dots, f_{n})$ is $\displaystyle \left(\sum_{j=0}^{n} f_j\right)^{n-1} f_0.$
\end{prop}

One sees that if $f_0 = f_1 = \cdots = f_{n} = 1$, $\MR_S(f_0, f_1, \dots, f_{r-1})$ is in bijection with $\R_S,$ and Proposition \ref{MRcard} recovers the result \eqref{RTcard}.
One famous way to prove \eqref{RTcard} is to construct the {\it Pr\"ufer sequence}. In fact, we will use this idea to prove Proposition \ref{MRcard}. Therefore, we will first review the construction of Pr\"ufer sequences.

\subsection*{Pr\"ufer sequences} Given a rooted tree $T \in \R_S,$ we define a sequence $T_1, T_2, \dots,$ $T_{n+1}$ of subtrees of $T$ as follows: Set $T_1 = T$. If $i<n+1$ and $T_i$ has been
defined, then define $T_{i+1}$ to be the tree obtained from $T_{i}$ by removing its largest leaf $v_i$ and the edge $e_i$ incident to $v_i$. Then define $w_i$ to be the other end of $e_i$, (i.e. $w_i$ is the parent of $v_i$), and let $\gamma(T):=(w_1,w_2, \ldots, w_{n})$. We call $\gamma(T)$ the {\it Pr\"ufer sequence} of $T.$

It is clear that $w_i \in S \cup \{0\}$ for $1 \leq i \leq n-1$ and $w_{n}=0$. Hence, $\gamma(T) \in (S \cup \{0\})^{n-1} \times \{0\}.$ The proof of the fact that $\gamma$ is a bijection from $\R_S$ to $(S \cup \{0\})^{n-1} \times \{0\}$  can be found in many places in the literature, for example, see \cite[Page 25]{stanleyec2}.

\begin{ex}
Let $T$ be the first tree shown in Figure \ref{fig_prufer}. Then $T_1,T_2,T_3$ and $T_4$ in Figure \ref{fig_prufer} are the first four trees appearing in the construction of the Pr\"{u}fer sequence of $T$. Continuing this construction, we obtain $\gamma(T)=(s_3,s_9,s_2,s_9,s_3,0,s_9,s_5,0)$.

\begin{figure}[h,t]
\begin{center}
\setlength{\unitlength}{8pt} \thicklines

\begin{picture}(44,11)
\put(0,0){
\begin{picture}(14,8)
\put(3.5,0){$T_1=T$}
\put(5,9){\circle*{0.3}}\put(5,9.2){$0$}
\put(5,9){\line(-3,-2){3}}
\put(5,9){\line(3,-2){3}}

\put(2,7){\circle*{0.3}}\put(0.8,7){$s_3$}
\put(2,7){\line(-3,-4){1.5}}
\put(2,7){\line(3,-4){1.5}}

\put(0.5,5){\circle*{0.3}}\put(-0.8,4.6){$s_8$} \put(0.5,5){\circle{0.8}}
\put(3.5,5){\circle*{0.3}}\put(3.8,4.6){$s_2$}
\put(3.5,5){\line(0,-1){2}}
\put(3.5,3){\circle*{0.3}}\put(3.6,2.2){$s_6$}

\multiput(8,7)(0,-2){3}{\circle*{0.3}}
\put(8.3,7){$s_5$}
\put(8.3,5){$s_9$}
\put(8,2){$s_1$}\put(6.3,2){$s_4$}
\put(9.3,2){$s_7$}
\put(8,7){\line(0,-1){4}}
\put(8,5){\line(-3,-4){1.5}}\put(6.5,3){\circle*{0.3}}
\put(8,5){\line(3,-4){1.5}}\put(9.5,3){\circle*{0.3}}
\end{picture}}

\put(11,0){
\begin{picture}(12,8)
\put(5,0){$T_2$}
\put(5,9){\circle*{0.3}}\put(5,9.2){$0$}
\put(5,9){\line(-3,-2){3}}
\put(5,9){\line(3,-2){3}}

\put(2,7){\circle*{0.3}}\put(0.8,7){$s_3$}
\put(2,7){\line(3,-4){1.5}}

\put(3.5,5){\circle*{0.3}}\put(3.8,4.6){$s_2$}
\put(3.5,5){\line(0,-1){2}}
\put(3.5,3){\circle*{0.3}}\put(3.6,2.2){$s_6$}

\multiput(8,7)(0,-2){3}{\circle*{0.3}}
\put(8.3,7){$s_5$}
\put(8.3,5){$s_9$}
\put(8,2){$s_1$}\put(6.3,2){$s_4$}
\put(9.3,2){$s_7$}\put(9.5,3){\circle{0.8}}
\put(8,7){\line(0,-1){4}}
\put(8,5){\line(-3,-4){1.5}}\put(6.5,3){\circle*{0.3}}
\put(8,5){\line(3,-4){1.5}}\put(9.5,3){\circle*{0.3}}
\end{picture}}

\put(24,0){
\begin{picture}(12,8)
\put(5,0){$T_3$}
\put(5,9){\circle*{0.3}}\put(5,9.2){$0$}
\put(5,9){\line(-3,-2){3}}
\put(5,9){\line(3,-2){3}}

\put(2,7){\circle*{0.3}}\put(0.8,7){$s_3$}
\put(2,7){\line(3,-4){1.5}}

\put(3.5,5){\circle*{0.3}}\put(3.8,4.6){$s_2$}
\put(3.5,5){\line(0,-1){2}}
\put(3.5,3){\circle*{0.3}}\put(3.6,2.2){$s_6$}\put(3.5,3){\circle{0.8}}

\multiput(8,7)(0,-2){3}{\circle*{0.3}}
\put(8.3,7){$s_5$}
\put(8.3,5){$s_9$}
\put(8,2){$s_1$}\put(6.3,2){$s_4$}
\put(8,7){\line(0,-1){4}}
\put(8,5){\line(-3,-4){1.5}}\put(6.5,3){\circle*{0.3}}
\end{picture}}

\put(35,0){
\begin{picture}(12,8)
\put(5,0){$T_4$}
\put(5,9){\circle*{0.3}}\put(5,9.2){$0$}
\put(5,9){\line(-3,-2){3}}
\put(5,9){\line(3,-2){3}}

\put(2,7){\circle*{0.3}}\put(0.8,7){$s_3$}
\put(2,7){\line(3,-4){1.5}}

\put(3.5,5){\circle*{0.3}}\put(3.8,4.6){$s_2$}

\multiput(8,7)(0,-2){3}{\circle*{0.3}}
\put(8.3,7){$s_5$}
\put(8.3,5){$s_9$}
\put(8,2){$s_1$}\put(6.3,2){$s_4$}
\put(8,7){\line(0,-1){4}}
\put(8,5){\line(-3,-4){1.5}}\put(6.5,3){\circle*{0.3}}\put(6.5,3){\circle{0.8}}
\end{picture}}
\end{picture}
\caption{Constructing the Pr\"{u}fer sequence of a rooted tree.}\label{fig_prufer}
\end{center}
\end{figure}
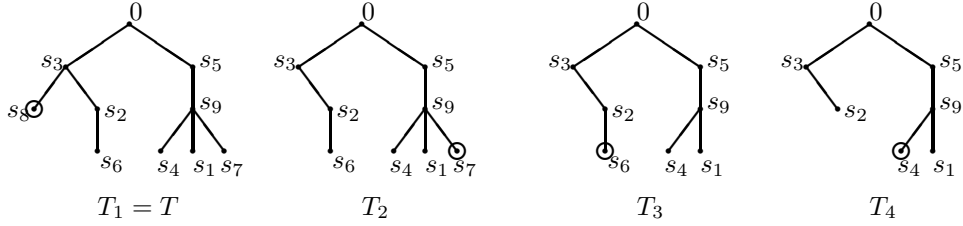

\end{ex}

\begin{proof}[Proof of Proposition \ref{MRcard}]
For convenience, we write $s_0 := 0.$ We denote by $H$ the set of matrices
$\displaystyle \left(\begin{array}{ccccc}w_1 & w_2 & \cdots & w_n \\ b_1 & b_2 & \cdots & b_n \end{array}\right)$
satisfying $(w_i, b_i) \in \bigcup_{j=0}^n \{ (s_j, k) \ | \ 1 \le k \le f_j\}$ for any $1 \le i \le n-1$ and $(w_n, b_n) \in \{ (0, k) \ | \ 1 \le k \le f_0\}.$ Since $\bigcup_{j=0}^n \{ (s_j, k) \ | \ 1 \le k \le f_j\}$ has cardinality $\sum_{j=0}^n f_j$ and $\{ (0, k) \ | \ 1 \le k \le f_0\}$ has cardinality $f_0,$ the cardinality of $H$ is $\left(\sum_{j=0}^{n} f_j\right)^{n-1} f_0.$
Our goal is to show that there is a bijection between $\MR_S(f_0, f_1, \dots, f_n)$ and $H.$

We will use the above algorithm for obtaining Pr\"ufer sequences of rooted trees to define this bijection.

Suppose $M = (T, \beta) \in \MR_S(f_0, f_1, \dots, f_n).$ Let $\gamma(T) = (w_1, \dots, w_n)$ be the Pr\"ufer sequence of $T$ and $e_1, \dots, e_n$ the edges removed in the procedure. We set $b_i := \beta(e_i).$ (In the labeled-edge representation of $G,$ $b_i$ is the label of the edge $e_i$ that is removed at step $i.$) Let
$$\tilde{\gamma}(M) = \left(\begin{array}{ccccc}w_1 & w_2 & \cdots & w_n \\ b_1 & b_2 & \cdots & b_n \end{array}\right).$$
One sees that $\tilde\gamma(G) \in H.$ Hence, $\tilde\gamma$ is a map from $\MR_S(f_0, f_1, \dots, f_n)$ to $H$.

On the other hand, suppose $\left(\begin{array}{ccccc}w_1 & w_2 & \cdots & w_n \\ b_1 & b_2 & \cdots & b_n \end{array}\right)$ is $H.$
Then $(w_1, \dots, w_n) \in (S \cup \{0\})^{n-1} \times \{0\}.$ Since $\gamma$ gives a bijection between $\R_S$ and $(S \cup \{0\})^{n-1} \times \{0\},$ we have $\gamma^{-1}(w_1, \dots, w_n) \in \R_S.$ Let $T := \gamma^{-1}(w_1, \dots, w_n).$ We can apply the algorithm to obtain the Pr\"ufer sequence of $T$ and record the order of the edges that were deleted. We then label the edge that was removed in the $i$th step with number $b_i.$ This procedure gives us a rooted tree $T$ with labeled edges, which is the edge-labeled representation of a multi-noded rooted tree $M = (T, \beta).$ One can check this procedure gives us the inverse of $\tilde\gamma.$

Therefore, $\tilde\gamma$ is a bijection between $\MR_S(f_0, f_1, \dots, f_n)$ and $H.$
Thus, the conclusion follows.
\end{proof}

\begin{ex}
For the multi-noded rooted tree $M$ in Example \ref{ex_multi-noded}, we have
$$\tilde{\gamma}(G) = \left(\begin{array}{ccccccccc}s_3 & s_9 & s_2 & s_9 & s_3 & 0 & s_9 & s_5 & 0 \\ 1 & 3 & 2 & 1 & 1 & 1 & 3 & 1 & 1  \end{array}\right).$$
\end{ex}

\begin{cor}\label{cor:card}
Suppose $\sum_{j=1}^{r-1} (e_j-1) = d-1.$
Then the cardinality of $\MR_S(1, e_1-1, \dots, e_{r-1}-1)$ is $d^{r-2}.$
\end{cor}
\begin{proof}
Let $n := r-1$ and $f_i := e_i -1$ for any $1 \le i \le n=r-1.$ Then by Proposition \ref{MRcard}, we have
\begin{equation}
|\MR_S(1, e_1-1, \dots, e_{r-1}-1)|=|\MR_S(1, f_1, \dots, f_n)| = (1 + \sum_{i=1}^n f_i)^{n-1} = d^{r-2}.
\end{equation}
\end{proof}

Corollary \ref{cor:card} provides us with a class of objects with cardinality $d^{r-2},$ which is the cardinality arising in Theorem \ref{main2}.
In the next two sections, we will describe in two steps a bijection between multi-noded rooted trees of vertex data $(1, e_1-1, \dots, e_{r-1}-1)$ and factorizations of a $d$-cycle of type $(e_1, \ldots, e_{r-1})$.

We will need the following definition when we construct the bijection in Section \ref{sec:g2t}.
\begin{defn}\label{defn:LMR}
Suppose $f_0, \dots, f_n$ are positive integers and let $d = \sum_{i=0}^n f_i.$ We say $(M, \l)$ is a {\it labeled multi-noded rooted tree of vertex data $(f_0,f_1,\dots,f_n)$} if $M \in \MR_S(f_0,f_1,\dots,f_n)$ and $\l$ is a labeling of the nodes of $M$ with set $[d].$ (So $\l$ is a bijection from the set of the nodes of $M$ to $[d].$)

We denote by $\LMR(f_0,f_1,\dots,f_n)$ the set of all the labeled multi-noded rooted trees of vertex data $(f_0,f_1,\dots,f_n).$
\end{defn}

The multi-noded rooted trees we define have connections to many other combinatorial objects. We finish this section with such an example. (This example is not related to the rest of the paper. The readers should feel free to skip it.)

\begin{cor}\label{cor:span}
Let $K_{m,n}$ be the complete bipartite graphs with vertex set $\{u_1, u_2, \ldots, u_m\}$ $\cup$ $\{v_1, v_2, \ldots, v_n\}$. Suppose $\ell_1, \ell_2, \ldots, \ell_m$ are positive integers with $\sum_{i=1}^{m} \ell_i=m+n-1$, then the number of spanning trees of $K_{m,n}$ such that $(u_1, u_2, \ldots, u_m)$ has degree sequence $(\ell_1, \ell_2, \ldots, \ell_m)$ is
\begin{equation}\label{spandeg}
n^{m-1}{n-1 \choose \ell_1 -1, \ell_2 -1, \ldots, \ell_m -1}
\end{equation}
and the total number of spanning trees of $K_{m,n}$ is
\begin{equation}\label{span}
n^{m-1}m^{n-1}.
\end{equation}
\end{cor}

The above enumeration result appears in \cite[Page 82, Ex 5.30]{stanleyec2}. We will provide a proof of it at the end of Section \ref{sec:g2t} using a bijection $\Phi^\L$ defined in that section. Here we just give an example of applying this result.

\begin{ex}
Suppose $m=2$ and $n=3$. There are $3$ spanning trees of $K_{2,3}$ with degree sequence $(u_1,u_2)=(1,3)$, $6$ spanning trees of $K_{2,3}$ with degree sequence $(u_1,u_2)=(2,2)$, and $3$ spanning trees of $K_{2,3}$ with degree sequence $(u_1,u_2)=(3,1)$. The total number of spanning threes of $K_{2,3}$ is $12=3^{2-1}2^{3-1}$.
\end{ex}

\section{Graphs associated to factorizations}\label{sec:bigraph}

Let $\tau \in S_d$ be a $d$-cycle, and $e_1, \dots, e_{r-1}$ integers no less than $2.$ Let $S = \{s_1<s_2< \cdots < s_{r-1}\}$ be a set of integers disjoint from $\{0,1,2,\dots,d\}.$
For any cycle $\gamma \in S_d,$ we denote by $\C_\gamma$ the circle with nodes labeled by numbers in $\gamma$ in clockwise order.

In this section, we associate a bipartite graph to each factorization of $\tau$ of type $(e_1,\dots,e_{r-1}).$ By discussing some properties of these graphs, we show that with the restriction $\sum_{j=1} (e_j-1) = d-1$ this association is an injection from $\Fac(d,r,\tau,e_1,\dots,e_{r-1})$ to its image set and thus is a bijection.

\begin{defn}
We call a graph $G$ an {\it $S$-$[d]$ bipartite graph} if the vertex set of $G$ is $S \cup [d]$ and any edge of $G$ connects a vertex in $S$ to a vertex in $[d].$ For any vertex $v$ in an $S$-$[d]$ bipartite graph $G$ or any subgraph of $G$, we call it an {\it $S$-vertex} if it is in $S$ and a {\it $[d]$-vertex} otherwise.

We denote by $\cG_S(d,r; e_1, \dots, e_{r-1})$ the set of all $S$-$[d]$ bipartite graphs $G$ satisfying for each $j: 1 \le j \le r-1$ the vertex $s_j$ has degree $e_j.$

\end{defn}

\subsection*{$S$-$[d]$ bipartite Graph associated to factorizations}
Suppose \eqref{req} and $(\sigma_1,$ $\dots,$ $\sigma_{r-1})$ is a factorization of $\tau$ of type $(e_1, \dots, e_{r-1})$. We associate to $(\sigma_1,$ $\dots,$ $\sigma_{r-1})$ a graph $G=(V,E)$ with vertex set $V = S \cup \supp(\tau) = S \cup [d]$ and edge set $E$ consisting of all the pairs $\{s_j, \nu\}$ where $\nu \in \supp(\sigma_j).$
We call $G$ a {\it factorization graph} of type $(d,r, \tau; e_1, \dots, e_{r-1})$.

\begin{ex}\label{ex_hurwitz}
Let $d=20$, $r=10$, $\tau=(1~2~\cdots~20)$ and $\sigma_1=(10~11)$, $\sigma_2=(14~15~19)$, $\sigma_3=(1~19)$, $\sigma_4=(3~4~5)$, $\sigma_5=(1~2~13)$, $\sigma_6=(15~16~17~18)$, $\sigma_7=(7~8~9~11)$, $\sigma_8=(19~20)$ and $\sigma_9=(2~5~6~11~12)$. One verifies that $(\sigma_1, \sigma_2, \cdots, \sigma_9)$ is a factorization of $\tau$ of type $(2,3,2,3,3,4,4,2,5)$. The corresponding factorization graph is shown in Figure \ref{fig_BiGraph}. (Note that the bipartite graph in the figure is drawn in a special way such that the $[d]$-vertices are embedded onto $\C_\tau.$ It will become clear later why we draw the graph this way.)

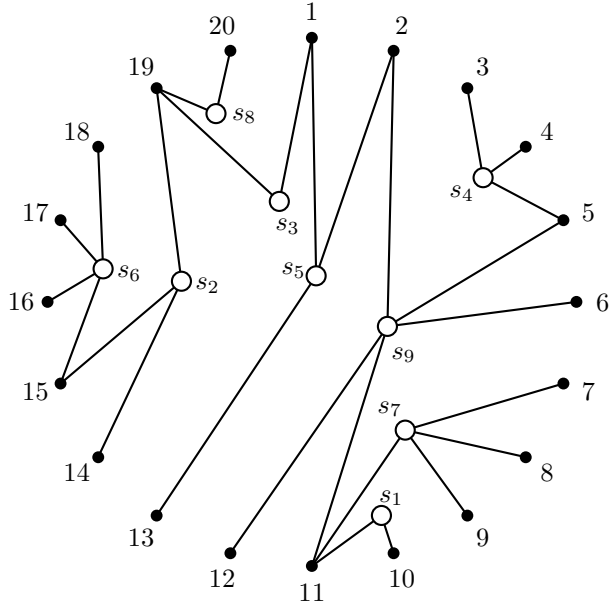
\begin{figure}[h,t] 
\begin{center}
\begin{pspicture}(-110pt,-110pt)(110pt,110pt)
\psset{unit=10pt}
\SpecialCoor
\degrees[360]


\rput(8;292){$s_1$}
\psline(8.5;288)(10;270) \psline(8.5;288)(10;288)
\pscircle[fillcolor=white,fillstyle=solid](8.5;288){4pt}

\rput(4;171){$s_2$} 
\psline(5;171)(10;126) \psline(5;171)(10;198) \psline(5;171)(10;216)
\pscircle[fillcolor=white,fillstyle=solid](5;171){4pt}

\rput(3;108){$s_3$} 
\psline(4;108)(10;90) \psline(4;108)(10;126)
\pscircle[fillcolor=white,fillstyle=solid](4;108){4pt}

\rput(7;36){$s_4$}
\psline(8;36)(10;18) \psline(8;36)(10;36) \psline(8;36)(10;54)
\pscircle[fillcolor=white,fillstyle=solid](8;36){4pt}

\rput(1.4;120){$s_5$}
\psline(1;81)(10;72) \psline(1;81)(10;90) \psline(1;81)(10;234)
\pscircle[fillcolor=white,fillstyle=solid](1;81){4pt}

\rput(7;171){$s_6$}
\psline(8;171)(10;144) \psline(8;171)(10;162) \psline(8;171)(10;180) \psline(8;171)(10;198)
\pscircle[fillcolor=white,fillstyle=solid](8;171){4pt}

\rput(5;306){$s_7$}
\psline(6;306)(10;270) \psline(6;306)(10;306) \psline(6;306)(10;324) \psline(6;306)(10;342)
\pscircle[fillcolor=white,fillstyle=solid](6;306){4pt}

\rput(7.5;110){$s_8$}
\psline(8;117)(10;126) \psline(8;117)(10;108)
\pscircle[fillcolor=white,fillstyle=solid](8;117){4pt}

\rput(4;330){$s_9$}
\psline(3;342)(10;0) \psline(3;342)(10;18) \psline(3;342)(10;72) \psline(3;342)(10;252) \psline(3;342)(10;270)
\pscircle[fillcolor=white,fillstyle=solid](3;342){4pt}

\psdots[linewidth=1.25pt](10;0) \rput(11;0){6}
\psdots[linewidth=1.25pt](10;18) \rput(11;18){5}
\psdots[linewidth=1.25pt](10;36) \rput(11;36){4}
\psdots[linewidth=1.25pt](10;54) \rput(11;54){3}
\psdots[linewidth=1.25pt](10;72) \rput(11;72){2}
\psdots[linewidth=1.25pt](10;90) \rput(11;90){1}
\psdots[linewidth=1.25pt](10;108) \rput(11;108){20}
\psdots[linewidth=1.25pt](10;126) \rput(11;126){19}
\psdots[linewidth=1.25pt](10;144) \rput(11;144){18}
\psdots[linewidth=1.25pt](10;162) \rput(11;162){17}
\psdots[linewidth=1.25pt](10;180) \rput(11;180){16}
\psdots[linewidth=1.25pt](10;198) \rput(11;198){15}
\psdots[linewidth=1.25pt](10;216) \rput(11;216){14}
\psdots[linewidth=1.25pt](10;234) \rput(11;234){13}
\psdots[linewidth=1.25pt](10;252) \rput(11;252){12}
\psdots[linewidth=1.25pt](10;270) \rput(11;270){11}
\psdots[linewidth=1.25pt](10;288) \rput(11;288){10}
\psdots[linewidth=1.25pt](10;306) \rput(11;306){9}
\psdots[linewidth=1.25pt](10;324) \rput(11;324){8}
\psdots[linewidth=1.25pt](10;342) \rput(11;342){7}
\end{pspicture}
\end{center}
\caption{A factorization graph $G$ of type $(20,10,(1~2~\cdots~20);$ $2,3,2,3,3,4,4,2,5)$.}\label{fig_BiGraph}
\end{figure}
\end{ex}

We denote by $\cG_S^*(d,r, \tau; e_1, \dots, e_{r-1})$ the set of all the factorization graphs of type $(d,r, \tau; e_1, \dots, e_{r-1})$.
Clearly $$\cG_S^*(d,r, \tau; e_1, \dots, e_{r-1}) \subset \cG_S(d,r; e_1, \dots, e_{r-1}).$$

One may notice that the factorization graph in Figure \ref{fig_BiGraph} is a tree. In fact this is not a coincidence. The following lemma and corollary discuss conditions when $G$ is a tree.

\begin{lem}\label{lem:tree}
Suppose $G \in \cG_S(d,r; e_1, \dots, e_{r-1})$ is connected. Then $\sum_{j=1}^{r-1} (e_j-1) = d-1$ if and only if $G$ is a tree.
\end{lem}
\begin{proof}
Any graph is a tree if and only if the graph is connected and the number of vertices is one more than the number of edges. Therefore, $G$ is a tree if and only if $1=|S\cup \supp(\tau)| - \sum_{j=1}^{r-1} e_j=r-1+d-\sum_{j=1}^{r-1} e_j,$ which is equivalent to $\sum_{j=1}^{r-1} (e_j-1) = d-1.$
\end{proof}

\begin{cor}\label{cor:tree}
Suppose $G \in \cG_S^*(d,r, \tau; e_1, \dots, e_{r-1}).$ Then $\sum_{j=1}^{r-1} (e_j-1) = d-1$ if and only if $G$ is a tree.
\end{cor}
\begin{proof}
Suppose $G$ is the factorization graph associated to $(\sigma_1, \dots, \sigma_{r-1})$, a factorization of $\tau.$ Since $\tau$ is a $d$-cycle, one sees $(\sigma_1, \dots, \sigma_{r-1})$ generates a transitive subgroup of $S_d.$ Thus, any two $[d]$-vertices of $G$ are connected by a path. However, any $S$-vertex is connected to some $[d]$-vertex. Hence, $G$ is connected. Then the conclusion follows from Lemma \ref{lem:tree}.
\end{proof}

It turns out that the same factorization graph can occur for different $\tau.$ However, it is true that any two different factorizations of a fixed $\tau$ have different factorization graphs, which is not obvious from the definition. We will show this fact at the end of this section by induction on $r$, and to achieve this, we discuss conditions on factorizations of $\tau$.

\begin{lem}\label{char2prd}
Suppose $\mu = (u_1, \dots, u_q)$ is a $q$-cycle and $\eta \in S_{\supp(\mu)}$ satisfying $\supp(\eta) = \{u_{j_1}, \dots, u_{j_p} \} \subseteq \{u_1, u_2, \ldots, u_q\}$ for some $j_1 > \cdots > j_p$. Let $s$ be the number of disjoint cycles (including the ones of length $1$) in the cycle decomposition of $\mu \eta$. Then $s \le p$, and the followings are equivalent:
\begin{ilist}
\itm $s = p.$
\itm $\eta = (u_{j_1}, \dots, u_{j_p})$.
\itm \label{2cycleprd} $\mu \eta = (u_{j_1+1}, u_{j_1+2}, \dots, u_q,  u_1, \dots,  u_{j_p}) (u_{j_p+1}, u_{j_p + 2}, \dots, u_{j_{p-1}}) \cdots$ \\
\indent \hspace{4mm} $\cdots (u_{j_2+1},  u_{j_2+2}, \dots,  u_{j_1}).$
\end{ilist}
\end{lem}

\begin{rem}
In this paper, whenever we talk about cycle decomposition, in addition to the disjoint cycles of length greater than $1$ appearing in the standard cycle decomposition of a permutation, we also include ``cycles'' of length $1.$ By convention, each of these contains exactly one fixed point of the permutation. We consider the support of each ``$1$-cycle'' to be its associated fixed point. Thus, the support of the cycles in the cycle decomposition of a permutation in $S_d$ always gives a partition of $[d].$
\end{rem}

\begin{proof}[Proof of Lemma \ref{char2prd}]

Clearly, if $u_i \not\in \supp(\eta)$, then $\mu(u_i) = \mu \eta(u_i).$ Hence, under the permutation $\mu \eta,$ we must have
\begin{eqnarray*}
& & u_{j_1+1} \mapsto  u_{j_1+2} \mapsto \cdots \mapsto u_q \mapsto  u_1 \mapsto u_2  \mapsto \cdots \mapsto  u_{j_p-1} \mapsto u_{j_p}\\
& & u_{j_p+1} \mapsto u_{j_p + 2} \mapsto \cdots \mapsto u_{j_{p-1}}, \\
& &u_{j_{p-1}+1}\mapsto u_{j_{p-1}+2}  \mapsto \cdots \mapsto u_{j_{p-2}}, \\
& & \cdots \\
& &u_{j_2+1} \mapsto u_{j_2+2} \mapsto  \cdots \mapsto u_{j_1}.
\end{eqnarray*}
Hence, the numbers in each line have to be in the same cycle in the cycle decomposition of $\mu \eta.$ Therefore the number of disjoint cycles in $\mu \eta$ is at most the number of lines we have above, i.e., $s \le p.$

It is easy to verify that (ii) and (iii) are equivalent. We show that (i) is equivalent to (iii).
We have $s = p$ if and only if the number at the end of each line is mapped to the number at the front under $\mu \eta.$ This means
$$\mu \eta(u_{j_{p}}) = u_{j_{1}+1}, \  \mu \eta(u_{j_{p-1}}) = u_{j_{p}+1}, \ \dots, \ \mu \eta(u_{j_1})= u_{j_2+1},$$
i.e.,
$$\eta(u_{j_{p}}) = \mu^{-1}(u_{j_{1}+1}) = u_{j_1},  \ \eta(u_{j_{p-1}}) = \mu^{-1}(u_{j_{p}+1}) = u_{j_p}, \  \dots $$
$$\dots, \ \eta(u_{j_1})= \mu^{-1}(u_{j_2+1}) = u_{j_2}.$$
Then our conclusion follows.
\end{proof}

\begin{rem}\label{rem:prd}
We can also understand Lemma \ref{char2prd} combinatorially: Suppose $\mu = (u_1, \dots, u_q)$ is a $q$-cycle and $\eta \in S_{\supp(\mu)}$ satisfying $\supp(\eta) = \{u_{j_1}, \dots, u_{j_p} \} \subseteq$ $\{u_1,$ $u_2,$ $\ldots,$ $u_q\}$ and $j_1 > \cdots > j_p$. Recall $\C_\mu$ is a circle whose nodes are labeled by $u_1,\dots,u_q$ in clockwise order.

Then the followings are equivalent:
\begin{ilist}
\itm There are $p$ cycles in the cycle decomposition of $\mu \eta.$
\itm The numbers in $\eta$ appear counterclockwise on $\C_\mu.$
\itm We can cut $\C_\mu$ into consecutive pieces such that each piece forms a cycle in the cycle decomposition of $\mu \eta$ when reading clockwise.
\end{ilist}
\end{rem}

\begin{ex}
Let $\mu = \tau=(1~2~\cdots ~20)$ and $\eta = \sigma_9^{-1}=(12~11~6~5~2)$ as in Example \ref{ex_hurwitz}. We have
\[\mu \eta = \tau \sigma_9^{-1} = (1~2~\cdots~20)(12~11~6~5~2)=(3~4~5)(7~8~9~10~11)(13~14~\cdots~20~1~2)(6)(12).\]
See Figure \ref{fig_circle}. Clearly $12,11,6,5,2$ appear counterclockwise on $\C_\tau.$
If we cut $\C_\tau$ after each of $2,5,6,11,12$, then we get exactly $5$  consecutive pieces $(3,4,5)$, $(6)$, $(7,8,9,10,11)$, $(12)$ and $(13,14, \ldots,20,1,2)$ when reading the numbers in clockwise order.

\begin{figure}[h,t] 
\begin{center}
\begin{pspicture}(-90pt,-90pt)(90pt,90pt)
\psset{unit=8pt}
\SpecialCoor
\degrees[360]

\pscircle(0;0){10}

\psdots[linewidth=1pt](10;0) \pscircle(10;0){4pt} \rput(11;0){6}
\psdots[linewidth=1pt](10;18) \pscircle(10;18){4pt} \rput(11;18){5}
\psdots[linewidth=1pt](10;36) \rput(11;36){4}
\psdots[linewidth=1pt](10;54) \rput(11;54){3}
\psdots*[linewidth=1pt](10;72) \pscircle(10;72){4pt} \rput(11;72){2}
\psdots[linewidth=1pt](10;90) \rput(11;90){1}
\psdots[linewidth=1pt](10;108) \rput(11;108){20}
\psdots[linewidth=1pt](10;126) \rput(11;126){19}
\psdots[linewidth=1pt](10;144) \rput(11;144){18}
\psdots[linewidth=1pt](10;162) \rput(11;162){17}
\psdots[linewidth=1pt](10;180) \rput(11.2;180){16}
\psdots[linewidth=1pt](10;198) \rput(11.2;198){15}
\psdots[linewidth=1pt](10;216) \rput(11.2;216){14}
\psdots[linewidth=1pt](10;234) \rput(11.2;234){13}
\psdots[linewidth=1pt](10;252) \pscircle(10;252){4pt} \rput(11.2;252){12}
\psdots[linewidth=1pt](10;270) \pscircle(10;270){4pt} \rput(11.2;270){11}
\psdots[linewidth=1pt](10;288) \rput(11;288){10}
\psdots[linewidth=1pt](10;306) \rput(11;306){9}
\psdots[linewidth=1pt](10;324) \rput(11;324){8}
\psdots[linewidth=1pt](10;342) \rput(11;342){7}

\psarc(0;0){8.8}{15}{60}\psline(9;14)(11;14)
\psarc(0;0){8.8}{67}{240}\psline(9;66)(11;66)
\psarc(0;0){8.8}{249}{258}\psline(9;248)(11;248)
\psarc(0;0){8.8}{267}{342}\psline(9;266)(11;266)
\psarc(0;0){8.8}{-3}{6}\psline(9;-4)(11;-4)
\end{pspicture}
\end{center}
\caption{The products of two cycles $\mu=(1~2~\cdots~20)$ and $\eta=(12~11~6~5~2)$.}\label{fig_circle}
\end{figure}
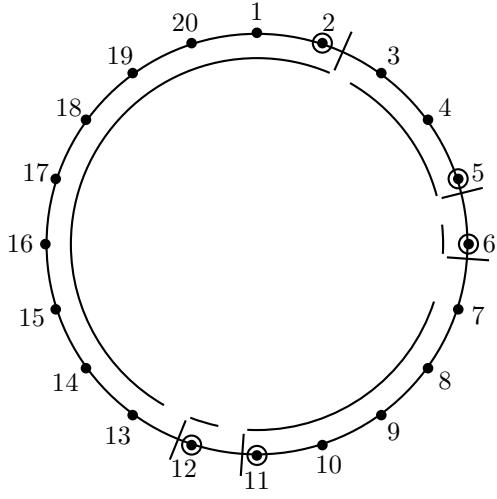
\end{ex}

\begin{lem}\label{chargraph}
  Suppose $\sum_{j=1}^{r-1} (e_j-1) = d-1$ and $(\sigma_1, \dots, \sigma_{r-1})$ is a factorization of $\tau$ of type $(e_1, \dots, e_{r-1}),$ and $G \in \cG_S^*(d,r, \tau; e_1, \dots, e_{r-1})$ is its corresponding factorization graph. Then $G$ is a tree by Corollary \ref{cor:tree}.

Suppose by deleting $s_{r-1}$ and its incident edges from $G,$ we obtain trees $Q_1,$ $\dots,$ $Q_k,$ $Q_{k+1}, \dots, Q_{e_{r-1}},$ where for $1 \le i \le k,$ the $[d]$-vertex set of $Q_i$ has size $m_i$ for some $m_i \ge 2,$ and for $k+1 \le i \le e_{r-1},$ $Q_i$ just contains one single $[d]$-vertex.

For any $i: 1 \le i \le k,$ let $B_i$ be the set of $j$ for which $s_j$ is in $Q_i.$ Then $\{B_1, \dots, B_k\}$ is a partition of $[r-2] :=\{1, 2,\dots, r-2\}.$

For $1 \le i \le k,$ let $\gamma_i := \prod_{j \in B_i} \sigma_j,$ where the product is taking over $j$ in increasing order, and for $k+1 \le i \le e_{r-1},$ let $\gamma_i$ be the $1$-cycle containing the only $[d]$-vertex of $Q_i.$  Then

\begin{ilist}
\itm $\gamma_1 \cdots \gamma_k \gamma_{k+1} \cdots \gamma_{e_{r-1}}$ is the cycle decomposition of $\sigma_1 \cdots \sigma_{r-2} = \tau \sigma_{r-1}^{-1}.$
\itm $\gamma_i$ is an $m_i$-cycle on the $[d]$-vertex set of $Q_i$.
\itm $(\sigma_j)_{j \in B_i}$ is a factorization of $\gamma_i,$ and $Q_i$ is the factorization graph associated to this factorization.
\itm $\sum_{j \in B_i} (e_j-1) = m_i -1.$
\end{ilist}

\end{lem}

\begin{ex}\label{ex_cycleprod}
Let $d,r,$ $\sigma_1,\dots,\sigma_9,\tau$ and $G$ be defined as in Example \ref{ex_hurwitz}. So $G$ is the graph in Figure \ref{fig_BiGraph}.
If we delete $s_9$ and its incident edges from $G$, we obtain $e_9=5$ trees, including two trees that are only a single $[d]$-vertex. Let $Q_1, Q_2$, $Q_3$, $Q_4$ and $Q_5$
denote the five trees with $[d]$-vertex set $\{3,4,5\}$, $\{7,8,\ldots,11\}$, $\{13,14,\ldots,20,1,2\}$, $\{6\}$ and $\{12\}$ respectively. Using the notation of Lemma \ref{chargraph},
we have $k=3$, $m_1=3$, $m_2=5$, $m_3=10$, and the corresponding partition of $[r-2]=[8]$ is $B_1=\{4\}$, $B_2=\{1,7\}$, $B_3=\{2,3,5,6,8\}$. Let $\gamma_1=\sigma_4=(3~4~5)$,
$\gamma_2=\sigma_1 \sigma_7=(7~8~\cdots~11)$, $\gamma_3=\sigma_2 \sigma_3 \sigma_5 \sigma_6 \sigma_8 =(13~14~\cdots~20~1~2)$, $\gamma_4=(6)$ and $\gamma_5=(12)$. One can check that
$\gamma_1\gamma_2\gamma_3\gamma_4\gamma_5$ is the cycle decomposition of $\tau\sigma_{9}^{-1}$, and for each $i: 1\leq i \leq 3$, (ii), (iii) and (iv) hold.

\end{ex}

\begin{proof}[Proof of Lemma \ref{chargraph}]
Since all the $e_j$'s are greater than $1,$ we have that any $Q_i$ for $k+1 \le i \le e_{r-1}$ does not contain any $S$-vertices. Therefore, each $s_j$ for any $j \in [r-2]$ is in one of $Q_1, \dots, Q_k.$ Thus, $\{B_1, \dots, B_k\}$ is a partition of $[r-2].$ Let $i \in \{1, \dots, k\}$ and $j \in B_i.$ One sees that all the $[d]$-vertices incident to $s_j$ have to be in $Q_i$ as well. Therefore, $\supp(\sigma_j)$ is contained in the $[d]$-vertex set of $Q_i.$
Thus, for any $j_1 \in B_{i_1}$ and $j_2 \in B_{i_2}$ with $i_1 \neq i_2$, we have that $\supp(\sigma_{j_1})$ and $\supp(\sigma_{j_2})$ are disjoint, which implies that $\sigma_{j_1} \sigma_{j_2} = \sigma_{j_2} \sigma_{j_1}.$ Hence,
$$\prod_{i=1}^k \gamma_i = \prod_{i=1}^k \prod_{j \in B_i} \sigma_j = \prod_{j=1}^{r-2} \sigma_j = \tau \sigma_{r-1}^{-1}.$$
Furthermore, for each $i: 1 \le i \le k,$ $\gamma_i = \prod_{j \in B_i} \sigma_j$ is a permutation on the $[d]$-vertex set of $Q_i.$
Therefore, the support of $\gamma_i$'s ($1 \le i \le e_{r-1}$) are completely disjoint. Hence, we can partition the cycles in the cycle decomposition of $\tau \sigma_{r-1}^{-1}$ into $e_{r-1}$ groups such that the product of the $i$th group of cycles is exactly $\gamma_i.$ This implies that $e_{r-1}$ is no greater than the number of cycles in the cycle decomposition of $\tau \sigma_{r-1}^{-1}.$ However, by applying Lemma \ref{char2prd} with $\mu = \tau$ and $\eta = \sigma_{r-1}^{-1},$ we have that the number of cycles in the cycle decomposition of $\tau \sigma_{r-1}^{-1}$ is no greater than $e_{r-1}.$ Hence, these two numbers are equal. So each $\gamma_i$ is one cycle in the cycle decomposition of $\tau \sigma_{r-1}^{-1}.$ We conclude (i),(ii) and (iii). Finally, (iv) follows from (iii) and Corollary \ref{cor:tree}.
\end{proof}

Combining Lemma \ref{chargraph}, Lemma \ref{char2prd} and Remark \ref{rem:prd}, we have the following corollary.

\begin{cor}\label{clockwise}
Suppose $\sum_{j=1}^{r-1} (e_j-1) = d-1$ and $(\sigma_1, \dots, \sigma_{r-1})$ is a factorization of $\tau$ of type $(e_1, \dots, e_{r-1}).$ Then we have the following conclusions:
\begin{ilist}
\itm For each $j: 1 \le j \le r-1,$ the numbers in $\sigma_j$ appear clockwise on $\C_\tau.$
\itm Let $\gamma_1, \dots, \gamma_k, \gamma_{k+1},\dots, \gamma_{e_{r-1}}$ be defined as in Lemma \ref{chargraph}. Then $\supp(\gamma_1),$ $\dots,$ $\supp(\gamma_{e_{r-1}})$ partition $\C_\tau$ into consecutive pieces. Furthermore, for each $i: 1 \le i \le e_{r-1},$ the numbers in $\gamma_i$ appear consecutively on $\C_\tau$ reading clockwise. Moreover, each $\gamma_i$ contains exactly one number from $\sigma_{r-1}$ and this number is the last number appearing on $\C_\tau$.
\end{ilist}
\end{cor}

\begin{proof}
By Lemma \ref{chargraph}, we have that the number of cycles in the cycle decomposition of $\tau \sigma_{r-1}^{-1}$ is equal to $e_{r-1},$ the size of the support of $\sigma_{r-1}^{-1}.$ Hence, by Lemma \ref{char2prd} and Remark \ref{rem:prd}, we have (ii) and the numbers in $\sigma_{r-1}$ appear clockwise on $\C_\tau$. We can conclude (i) for other $j$'s by applying Lemma \ref{chargraph}/(iii)(iv), Lemma \ref{char2prd} and Remark \ref{rem:prd} recursively.
\end{proof}

By Corollary \ref{clockwise}/(i), one sees that with the condition $\sum_{j=1}^{r-1} (e_j-1) =d-1,$ no two different factorizations of $\tau$ can have the same factorization graph.

\begin{cor}\label{bijFacGraph}
Suppose $\sum_{j=1}^{r-1} (e_j-1) =d-1.$ The way we associate a graph to a factorization gives a bijection between the set $\Fac(d,r,\tau;e_1, \dots, e_{r-1})$ and the set $\cG_S^*(d,r, \tau; e_1, \dots, e_{r-1})$.
\end{cor}

\begin{rem}
We remark that if $e_1 = \cdots = e_{r-1} = 2$, then $d = r$ and $\Fac(d,d,\tau;$ $2, \dots, 2)$ contains factorizations of a $d$-cycle $\tau$ into $d-1$ transpositions. In this case
for any $G \in \cG_S^*(d,d,\tau;2,2,\ldots,2)$, the $S$-vertices of $G$ have degree 2. For each $S$-vertex $s_j \in G$, suppose $s_j$ is incident to $\nu_{j_1}$ and $\nu_{j_2}$. We can replace $s_j$ and its two incident edges by one edge connecting $\nu_{j_1}$ and $\nu_{j_2}.$ Then we get a tree on vertex set $[d]$. Therefore, the bijection discussed in Corollary \ref{bijFacGraph} becomes a bijection between trees on $d$ vertices and factorizations of a $d$-cycle into $d-1$ transpositions, which is the same as the bijection defined by Moszkowski in \cite{Moszkowski} and the \emph{circle chord diagram} construction defined by Goulden and Yong in \cite{Treelike}.
\end{rem}

\section{A Bijection between factorization graphs and multi-noded rooted trees}\label{sec:g2t}

For convenience, we assume $\tau = (1~2~\cdots~d).$
In this section, we will define a map from factorization graphs in $\cG_S^*(d,r, \tau=(1~\cdots ~d); e_1, \dots, e_{r-1})$ to multi-noded rooted trees in $\MR_S(1, e_1-1, \dots, e_{r-1}-1)$, which we will show later is a bijection assuming $\sum_{j=1}^{r-1} (e_j-1) = d-1$. Clearly, such a bijection can be extended to any $\tau.$

We now construct our map. 

\begin{defn}Assume $\sum_{j=1}^{r-1} (e_j-1) = d-1$.
For any $G \in \cG_S^*(d,r, (1~2~\cdots ~d);$ $e_1, \dots, e_{r-1})$, we have that $G$ is a tree by Corollary \ref{cor:tree}. We make the $[d]$-vertex $1$ of $G$ a root, and call the resulting rooted tree $G^\R$.
It is clear that $s_i$ has $e_i-1$ children in $G^\R,$ for each $i: 1 \le i \le r-1.$

Recall that labeled multi-noded rooted trees are defined in Definition \ref{defn:LMR}. We define $\Phi^\L(G)=(M,{\mathfrak l})$ to be the labeled multi-noded rooted tree, where $M$ is in its multi-noded representation obtained from $G^\R$ in the following way:
\begin{alist}
\itm We make the root $1$ of $G^\R$ a single-noded vertex, which is the root of $\Phi^\L(G).$ We keep the node label $1$ and label the single-noded vertex with $s_0 = 0.$
\itm For each $i: 1 \le i \le r-1,$ suppose $\nu_1 < \dots < \nu_{e_i-1}$ are the children of $s_i$ and $\nu$ is the parent of $s_i$ in $G^\R.$  Let $s_i$ be an $(e_i-1)$-noded vertex containing nodes which are labeled by $\nu_1, \dots, \nu_{e_i-1}$ from left to right. Then connect $s_i$ to the node $\nu.$
\end{alist}
One sees that $\Phi^\L(G)=(M, \l)$ is in $\LMR_S(1, e_1-1, \dots, e_{r-1}-1),$ where $M \in \MR_S(1, e_1-1, \dots, e_{r-1}-1)$. We denote $M$ by $\Phi(G).$

\end{defn}

\begin{ex}
Let $d=20$, $r=10$, $\tau=(1~2~\cdots~20)$ and $G$ be the graph shown in Figure \ref{fig_BiGraph}, which is the bipartite graph associated to the factorization defined in Example \ref{ex_hurwitz}.
Then $G \in \cG_S^*(d,r,\tau;2,3,2,3,3,4,4,2,5)$ and $\Phi^\L(G)$ is
a labeled multi-noded rooted tree in $\LMR_S(1, 1, 2, 1, 2, 2, 3, 3, 1, 4)$. Figure \ref{fig_MutiLabelled} shows the multi-noded representation of $\Phi^\L(G).$
After removing labels for the nodes, we get $\Phi(G)$, which is the multi-noded rooted tree shown in Figure \ref{tworep}(b).

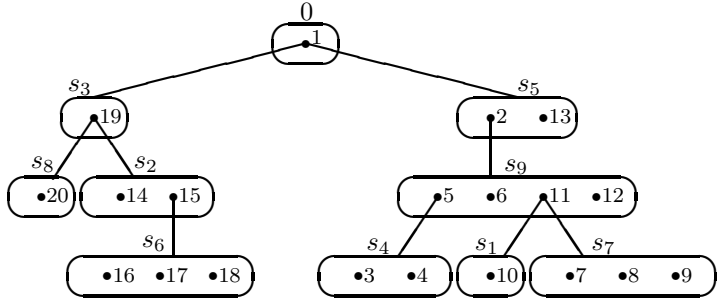
\begin{figure}[h,t]
\begin{center}
\setlength{\unitlength}{10pt} \thicklines
\begin{picture}(26,11)
\put(10,9.8){\circle*{0.3}}\put(9.8,10.7){$0$}
\put(10,9.8){\oval(2.5,1.5)}
\put(10,9.8){\line(-4,-1){8}}
\put(10,9.8){\line(4,-1){8}}
\put(10.2,9.8){\footnotesize{1}}

\put(2,7){\circle*{0.3}}\put(1,8){$s_3$}
\put(2,7){\line(-2,-3){1.55}}
\put(2,7){\line(2,-3){1.5}}
\put(2,7){\oval(2.5,1.5)}
\put(2.2,6.8){\footnotesize{19}}

\put(0,4){\circle*{0.3}}\put(-0.4,5.1){$s_8$}
\put(0,4){\oval(2.5,1.5)}
\put(0.2,3.8){\footnotesize{20}}

\multiput(3,4)(2,0){2}{\circle*{0.3}}
\put(3.5,5.1){$s_2$}
\put(4,4){\oval(5,1.5)}
\put(5,4){\line(0,-1){2.25}}
\put(3.2,3.8){\footnotesize{14}}
\put(5.2,3.8){\footnotesize{15}}

\multiput(2.5,1)(2,0){3}{\circle*{0.3}}
\put(3.8,2.1){$s_6$}
\put(4.5,1){\oval(7,1.5)}
\put(2.7,0.8){\footnotesize{16}}
\put(4.7,0.8){\footnotesize{17}}
\put(6.7,0.8){\footnotesize{18}}

\multiput(17,7)(2,0){2}{\circle*{0.3}}\put(18,8){$s_5$}
\put(18,7){\oval(4.5,1.5)}
\put(17,7){\line(0,-1){2.25}}
\put(15,4){\line(-2,-3){1.5}}
\put(17.2,6.8){\footnotesize{2}}
\put(19.2,6.8){\footnotesize{13}}

\multiput(15,4)(2,0){4}{\circle*{0.3}}\put(18,4){\oval(9,1.5)}
\put(17.5,5){$s_9$}
\put(15.2,3.8){\footnotesize{5}}
\put(17.2,3.8){\footnotesize{6}}
\put(19.2,3.8){\footnotesize{11}}
\put(21.2,3.8){\footnotesize{12}}

\multiput(12,1)(2,0){2}{\circle*{0.3}}\put(13,1){\oval(5,1.5)}
\put(12.3,2){$s_4$}
\put(12.2,0.8){\footnotesize{3}}
\put(14.2,0.8){\footnotesize{4}}

\put(19,4){\line(-2,-3){1.5}}
\put(19,4){\line(2,-3){1.5}}

\put(17,1){\circle*{0.3}}
\put(17,1){\oval(2.5,1.5)}
\put(16.4,2){$s_1$}
\put(17.2,0.8){\footnotesize{10}}

\multiput(20,1)(2,0){3}{\circle*{0.3}} \put(22,1){\oval(7,1.5)}
\put(20.8,2){$s_7$}
\put(20.2,0.8){\footnotesize{7}}
\put(22.2,0.8){\footnotesize{8}}
\put(24.2,0.8){\footnotesize{9}}

\end{picture}
\caption{A multi-noded rooted tree with labeled nodes.}\label{fig_MutiLabelled}
\end{center}
\end{figure}
\end{ex}

\begin{thm}\label{thm:bijg2t}
Suppose $\sum_{j=1}^{r-1} (e_j-1) = d-1$.
Then $\Phi$ gives a bijection from $\cG_S^*(d,r,(1~2~\cdots~d);e_1,\dots,e_{r-1})$ to $\MR_S(1, e_1-1, \dots, e_{r-1}-1).$
\end{thm}

Given this theorem, we can prove our main theorems.

\begin{proof}[Proof of Theorems \ref{main2} and \ref{main1}]
	By Corollary \ref{bijFacGraph} and Theorem \ref{thm:bijg2t}, we see that we obtain a (two-step) bijection from factorizations of $\tau=(1~2~\cdots~d)$ of type $(e_1, \dots, e_{r-1})$ to multi-noded rooted trees of vertex data $(1, e_1-1, \dots, e_{r-1}-1):$
	\begin{eqnarray*}
		\Fac(d,r,(1~2~\cdots~d);e_1,\dots,e_{r-1}) &\longrightarrow& \cG_S^*(d,r,(1~2~\cdots~d);e_1,\dots,e_{r-1}) \\
		&\longrightarrow& \MR_S(1, e_1-1, \dots, e_{r-1}-1).
	\end{eqnarray*}
	We can extend this bijection to any $\tau.$ Hence, Theorem \ref{main2} follows from Corollary \ref{cor:card}.
As we discussed in the introduction, Theorem \ref{main1} follows.
\end{proof}

We remark that although the bijection we construct has two steps, each of the steps is quite simple.

We finish this section with a proof of Corollary \ref{cor:span}.

\begin{proof} [Proof of Corollary \ref{cor:span}.]
For every such spanning tree, we can get a labeled multi-rooted tree $(M, \l)$ using a map similar to the bijection $\Phi^{\L}$ we defined for factorization graphs, except that we make $v_1$ the single-noded root instead of 1. Here $M$ has vertex data $(1, \ell_1-1, \ell_2 -1, \ldots, \ell_m -1)$, and $\l$ can be any labeling such that the single-noded root is labeled $v_1$ and all the other nodes are labeled with $u_1, u_2, \ldots, u_m$. There are $(1 + \sum_{i=1}^m (\ell_i -1))^{m-1}\cdot 1=n^{m-1}$ such $M$'s, and for each $M$, there are ${n-1 \choose \ell_1 -1, \ell_2 -1, \ldots, \ell_m -1}$ ways to label it, therefore we get \eqref{spandeg}. (Note that if $\ell_i=1$ for some $i=1, 2, \ldots, m$, we just neglect the ``empty" vertex and the result still holds.) By summing over all positive integer solutions of $\ell_1 +\ell_2+ \cdots + \ell_m=n+m-1$, we get \eqref{span}.
\end{proof}

We devote the rest of the paper to the proof of Theorem \ref{thm:bijg2t}.

\section{Characterization of factorization graphs}\label{sec:charGraph}

In this section, we will give a proposition (Proposition \ref{prop:bijgraph}) to characterize properties of graphs in $\cG_S^*(d,r, \tau; e_1, \dots, e_{r-1})$, which will be used to prove Theorem \ref{thm:bijg2t} in the next section.
We first give definitions that are useful for the statement of the proposition.

\begin{defn}
Suppose $S' \subseteq S$ and $\gamma$ is a cycle in $S_d.$
Let $G$ be an $S'$-$\supp(\gamma)$ bipartite tree.

Suppose $s \in S'.$ We say $s$ has the {\it consecutive partition property (or CPP) on $(G,\gamma)$} if after we remove $s$ and all its incident edges from $G$, the sets of $[d]$-vertices of the subtrees we obtain partition the circle $\C_\gamma$ into consecutive pieces.

Suppose $\nu \in \supp(\gamma)$ and $\{s_{j_1} < s_{j_2} < \cdots < s_{j_t}\}$ are the set of $S$-vertices incident to $\nu$ in $G$. By removing $\nu$ and all its incident edges, suppose we obtain $t$ subtrees. We say $\nu$ has the {\it counterclockwise increasing consecutive partition property (or CICPP) on $(G,\gamma)$} if the following are satisfied:
\begin{alist}
\itm The $[d]$-vertices of the $t$ subtrees partition $\C_\gamma \setminus \{\nu\}$ into consecutive pieces.
\itm If we order the pieces in counterclockwise order on $\C_\gamma$ starting from $\nu,$ then the $m$-th piece is the $[d]$-vertex set of the subtree that contains vertex $s_{j_m}$ for any $1 \le m \le t.$
\end{alist}
\end{defn}

We can restate part of Corollary \ref{clockwise}/(ii) with this definition using the connection between $\gamma_i$ and $Q_i$ discussed in Lemma \ref{chargraph}.
\begin{cor}\label{sr-1CPP}
Suppose $\sum_{j=1}^{r-1} (e_j-1) = d-1$ and $G \in \cG_S^*(d,r, \tau; e_1, \dots, e_{r-1}).$ Then $s_{r-1}$ has CPP on $(G, \tau).$
\end{cor}

The properties CPP and CICPP are not independent. In fact we have the following lemma.

\begin{lem}\label{CI=>CPP}
Suppose $S' \subseteq S$ and $\gamma$ is a cycle in $S_d$.
Let $G$ be an $S'$-$\supp(\gamma)$ bipartite tree. Suppose $s \in S'$. If all the $[d]$-vertiecs incident to $s$ have CICPP on $(G,\gamma)$, then $s$ has CPP on $(G,\gamma)$.
\end{lem}

\begin{proof}

Suppose $\nu_1,\dots,\nu_k$ are the $[d]$-vertices incident to $s.$ Let $Q_1, \dots, Q_k$ be the subtrees containing $\nu_1,\dots,\nu_k$ respectively obtained from $G$ by removing $s$ and its incident edges.
One sees that it suffices to show that for each $i: 1 \le i \le k,$ the union of $[d]$-vertex sets of $Q_{i'}$ with $i' \neq i$ is a consecutive piece on $\C_\gamma.$ However, this follows from that $\nu_i$ has CICPP on $(G, \gamma)$ since this union is exactly the $[d]$-vertex set of the tree containing $s$ obtained by deleting the edge $\{s, \nu_i\}$ from $G.$
\end{proof}

We now state the main result of this section.
\begin{prop}\label{prop:bijgraph}
Suppose $G \in \cG_S(d,r; e_1, \dots, e_{r-1})$.

Then $\sum_{j=1}^{r-1} (e_j-1) = d-1$ and $G \in \cG_S^*(d,r, \tau; e_1, \dots, e_{r-1})$ if and only if $G$ satisfies the following conditions:
\begin{enumerate}
\item $G$ is a tree.
\item Any $[d]$-vertex of $G$ has CICPP on $(G, \tau).$
\end{enumerate}
Therefore, by Lemma \ref{CI=>CPP}, if $\sum_{j=1}^{r-1} (e_j-1) = d-1$ and $G \in \cG_S^*(d,r, \tau; e_1, \dots, e_{r-1})$, we also have the following:
\begin{enumerate}
\item[(3)] Any $S$-vertex of $G$ has CPP on $(G, \tau)$.
\end{enumerate}
\end{prop}

One can check that (3) of Proposition \ref{prop:bijgraph} is equivalent to the condition that $G$ has a planar embedding with $[d]$-vertices on the circle $C_\tau$ and $S$-vertices inside the circle.

\begin{ex}
Let $G$ be the graph shown in Figure \ref{fig_BiGraph}, which is the bipartite graph associated to the factorization defined in Example \ref{ex_hurwitz}.

From Example \ref{ex_cycleprod}, we see that $s_{9}$ has CPP on $(G, (1~2~\cdots~20))$, where the corresponding partition is $\{\{3,4,5\}$,$\{6\}$,$\{7,8,\ldots,11\}$,$\{12\}$,$\{13,14,\ldots,20,$ $1,2\}\}$.

If we remove the $[d]$-vertex $19$ and all its incident edges, we get three trees $T_1$, $T_2$ and $T_3$ whose vertex sets are $\{s_2,s_6\}\cup\{14, 15, \ldots, 18\}$, $\{s_1,s_3,s_4,s_5,s_7,s_9\}\cup\{1,2,\ldots,13\}$ and $\{s_8\}\cup\{20\}$, respectively. It is easy to see that the $[d]$-vertex sets of $T_1$, $T_2$ and $T_3$ partition the circle $(1~2~\cdots~18~20)$ into consecutive pieces, and these pieces are in counterclockwise order on the circle starting from 19. Moreover, the $S$-vertices incident to 19 are $s_2$, $s_6$ and $s_8$, and satisfy that $s_2 \in T_1$ $s_6 \in T_2$ and $s_8 \in T_3$. Thus 19 has CICPP on $(G, (1~2~\cdots~20))$.

The readers can check that all the other $S$-vertices have CPP on $(G, (1~2~\cdots~20))$, and all the other $[d]$-vertices have CICPP on $(G, (1 ~2~\cdots ~20))$.
\end{ex}

We will use the following lemma to prove Proposition \ref{prop:bijgraph}.
\begin{lem}\label{lem:iffCPP}
Suppose $G$ is an $S$-$[d]$ bipartite tree. Let $\nu_0$ be a $[d]$-vertex of $G$.

Suppose $Q$ and $\bar{Q}$ are two subtrees of $G$ satisfying: (1) The union of $Q$ and $\bar{Q}$ is $G$; (2) $\nu_0$ is the only common vertex of $Q$ and $\bar{Q}$; (3) the $[d]$-vertex set of $Q$ is a consecutive piece on $\C_\tau$ and ends with $\nu_0$ when reading clockwise.

Let $\gamma$ be the cycle obtained by reading the $[d]$-vertices of $Q$ in clockwise order on $\C_\tau$. Then we have the following:
\begin{ilist}
\itm For any $\nu \neq \nu_0$ a $[d]$-vertex of $Q,$ $\nu$ has CICPP on $(Q,\gamma)$ if and only if $\nu$ has CICPP on $(G,\tau).$
\end{ilist}

If we suppose further that $\{s_{j_1} < s_{j_2} < \cdots < s_{j_t}\}$ are the set of $S$-vertices incident to $\nu_0$ in $G,$ and $s_{j_1},\dots, s_{j_{t-1}}$ are in $Q$ and $s_{j_t}$ is in $\bar{Q},$ then
\begin{ilist}
\item[(ii)] $\nu_0$ has CICPP on $(Q,\gamma)$ if and only if $\nu_0$ has CICPP on $(G,\tau).$
\end{ilist}

\end{lem}

\begin{proof}
The labeling of $\tau$ does not matter, so we can always relabel so that $\nu_0$ is the size of the $[d]$-vertex set of $Q.$ Without loss of generality, we assume $\tau = (1~2~\cdots~d)$ and the $[d]$-vertex set of $Q$ is $\{1,2,\dots,\nu_0\}.$ So the $[d]$-vertex set of $\bar{Q}$ is $\{\nu_0, \nu_0+1, \dots, d\}.$ We also let $\gamma$ be the cycle $(1~2~\cdots~\nu_0).$

\begin{ilist}

\itm Suppose by removing $\nu$ and its incident edges from $G$, we get trees $T_1, \dots,$ $T_{t}.$ We can assume $T_1$ is the tree that contains $\nu_0.$ Let $T_1'$ be the tree obtained from $T_1$ by deleting $\bar{Q}.$ One can check that $T_1', T_2, \dots, T_{t}$ are the trees we obtain by removing $\nu$ and its incident edges from $Q$.

Suppose $\nu$ has CICPP on $(Q,\gamma)$. Then the $[d]$-vertices of $T_1', T_2, \dots, T_{t}$ partition $\C_\gamma \backslash \{\nu\}$ into consecutive pieces. Because $T_1'$ contains $\nu_0,$ the $[d]$-vertex set of $T_1'$ is of the form $\{\alpha,\alpha+1, \cdots, \nu_0, 1, 2,\cdots,\beta\}$ for some $0 \le \beta < \alpha \le \nu_0,$ and the other $t-1$ trees partition $[\beta+1, \alpha-1]\setminus\{\nu\}$ into consecutive pieces. However, the $[d]$-vertices of $\bar{Q}$ are $\{\nu_0,\nu_0+1,\dots, d\}.$ Hence, the $[d]$-vertices of $T_1$ are $\{\alpha, \alpha+1, \dots, d, 1, 2,\dots, \beta\}$. Therefore, the $[d]$-vertices of $T_1, T_2, \dots, T_{t}$ partition $\C_\tau \backslash \{\nu\}$ into consecutive pieces.
Moreover,  condition b) in the definition of $\nu$ having CICPP on $(G,\tau)$ can also be verified.
Therefore we proved that $\nu$ has CICPP on $(G, \tau)$.

By similar arguments we can prove the other direction that if $\nu$ has CICPP on $(G,\tau),$ then $\nu$ has CICPP on $(Q,\gamma)$.

\itm Let $T_1, \dots, T_t$ be the subtrees obtained from $G$ by removing $\nu_0$ and its incident edges, where $T_m$ contains $s_{j_m}$ for each $m: 1 \le m \le t.$ One checks that $\bar{Q}$ is the union of $T_t$ and the edge $\{\nu_0, s_{j_t}\}$ and $Q$ is the union of $T_1, \dots, T_{t-1}$ and edges $\{\{\nu_0,s_{j_m}\}\}_{m=1}^{t-1}.$ Hence, $T_1, \dots, T_{t-1}$ are the trees we obtain by removing $\nu_0$ and its incident edges from $Q$, and the $[d]$-vertex set of $T_t$ is $\{\nu_0+1, \nu_0+2, \dots, d\}.$
Now it is easy to verify that $\nu_0$ has CICPP on $(Q,\gamma)$ if and only if $\nu_0$ has CICPP on $(G,\tau).$
\end{ilist}
\end{proof}

\begin{proof}[Proof of Proposition \ref{prop:bijgraph}]
We prove the proposition by induction on $r.$ Suppose $r=2.$ The condition $\sum_{j=1}^{r-1} (e_j-1) =d-1$ is equivalent to $e_1 = d.$ Under this condition, $\cG_S^*(d,r, \tau; e_1, \dots, e_{r-1})$ contains only one graph $G_0=(\{s_1\}\cup [d], \{\{s_1,\nu\}\}_{\nu=1}^d),$ which satisfies (1) and (2). On the other hand, if $G$ satisfies (1) and (2), one sees that $G=G_0$, which is in $\cG_S^*(d,r, \tau; e_1, \dots, e_{r-1})$. Furthermore, we have to have $d = e_1.$

Suppose $r_0 \ge 3$ and the proposition holds for any $r < r_0.$ We prove the case $r = r_0.$ Let $G=(V,E) \in \cG_S(d,r;e_1,\dots,e_{r-1}).$ For convenience, for each $j: 1 \le j \le r-1$, we define the following:
\begin{itemize}
\itm Let $E_j$ be the set of edges in $G$ that are incident to $s_j.$
\itm Let $P_{j}$ be the ``star-shaped'' graph whose vertices are $s_{j}$ and the $e_{j}$ $[d]$-vertices incident to $s_{j},$ and whose edge set is $E_{j}.$
\end{itemize}

Suppose $\sum_{j=1}^{r-1} (e_j-1) = d-1$ and $G$ is the graph associated to a factorization $(\sigma_1, \dots, \sigma_{r-1}).$ Then $G$ is a tree by Corollary \ref{cor:tree}.
We only need to show (2).
Let $k, Q_1, \dots, Q_k, Q_{k+1}, \dots, Q_{e_{r-1}},$ $\gamma_1, \dots, \gamma_k, \gamma_{k+1},\cdots, \gamma_{e_{r-1}}$ and $B_1, \dots, B_k$ be defined as in Lemma \ref{chargraph}. By Lemma \ref{chargraph}/(iii),(iv) and the induction hypothesis, we know that $Q_i$ satisfies (1) and (2) for $1 \le i \le k.$

For any $i: 1 \le i \le k,$ we define $\nu_i$ to be the $[d]$-vertex of $Q_i$ that was incident to $s_{r-1}$ and $\bar{Q}_i$ the union of $P_{r-1}$ and $\cup_{i'\neq i} Q_{i'}.$
One checks that the union of $Q_i$ and $\bar{Q}_i$ is $G$ and $\nu_i$ is the only common vertex of $Q_i$ and $\bar{Q}_i.$
Thus, using these together with Corollary \ref{clockwise}/(ii), one sees that the hypothesis for (i) of Lemma \ref{lem:iffCPP} is satisfied by setting $Q= Q_i$, $\bar{Q} = \bar{Q}_i$ and $\gamma = \gamma_i.$

Let $\nu$ be an $S$-vertex of $G.$ Suppose $\nu$ is not in $\supp(\sigma_{r-1}),$ the set of vertices incident to $s_{r-1}.$ Then $\nu$ is in $Q_i$ for some $i: 1 \le i \le k.$ Since $\nu$ has CICPP on $(Q_i,\gamma_i)$, by Lemma \ref{lem:iffCPP}/(i), $\nu$ has CICPP on $(G,\tau).$ Suppose $\nu$ is in $\supp(\sigma_{r-1}).$ Then $\nu \in Q_i$ for some $i: 1 \le i \le e_{r-1}.$ If $i > k,$ $\nu$ is the only vertex in $Q_i$ and $s_{r-1}$ is the only vertex that is incident to $\nu.$ Then $\nu$ automatically has CICPP on $(G,\tau).$ Suppose $i \le k.$ Since $s_{r-1}$ is the biggest $S$-vertex incident to $\nu$, the conclusion follows from Lemma \ref{lem:iffCPP}/(ii) and the fact that $\nu$ has CICPP on $(Q_i,\gamma_i).$

Therefore, we proved that if $\sum_{j=1}^{r-1} (e_j-1) = d-1$ and $G \in \cG_S^*(d,r,\tau,e_1, \dots, e_{r-1}),$ then $G$ satisfies (1) and (2).

Suppose $G$ satisfies (1) and (2). Since $G$ is a tree which is connected, by Lemma \ref{lem:tree}, $\sum_{j=1}^{r-1}(e_j-1)=d-1.$ Hence, we only need to prove that $G$ is a factorization graph of type $(d,r,\tau;e_1,\dots, e_{r-1}).$
For each $j: 1 \le j \le r-1,$ we define $\sigma_j$ to be the $e_j$-cycle obtained by reading $[d]$-vertices incident to $s_j$ in clockwise order as appeared in $\C_\tau$.
It suffices to show that $\sigma_1 \cdots \sigma_{r-1} = \tau.$

We assume the $[d]$-vertices incident to $s_{r-1}$ are $\nu_1,$ $\dots,$ $\nu_{e_{r-1}}.$ Let $Q_1, \dots, Q_{e_{r-1}}$ be the subtrees we obtain by deleting $s_{r-1}$ and its incident edges from $G$, where $Q_i$ contains $\nu_i$ for each $i.$ Since $s_{r-1}$ has CPP on ($G$,$\tau$), the $[d]$-vertex set of $Q_i$ is a consecutive piece on $\C_\tau$ containing $\nu_i.$ We claim that
\begin{ilist}
\itm the $[d]$-vertex set of $Q_i$ is a consecutive piece on $\C_\tau$ which ends with $\nu_i$ when read in clockwise order, for each $i: 1 \le i \le e_{r-1}.$
\end{ilist}
One sees that it is enough to prove that
\begin{ilist}
\item[(i$'$)] the $[d]$-vertex set of $Q_i$ does not contain $\tau(\nu_i),$ the number after $\nu_i$ on $\C_\tau$ in clockwise order, for each $i.$
\end{ilist}
We assume to the contrary that for some $i,$ the $[d]$-vertex set of $Q_i$ contains $\tau(\nu_i).$ Then among the subtrees we obtain by removing $\nu_i$ and its incidence edges, the one containing $s_{r-1}$ does not contain the vertex $\tau(\nu_i)$, which contradicts the assumption that $\nu_i$ has CICPP on $(G,\tau)$. Therefore, (i$'$) holds and thus (i) holds.

Let $m_i$ be the size of the $[d]$-vertex set of $Q_i$ for each $i.$ Without loss of generality, we may assume $m_1,\dots,m_k \ge 2$ and $m_{k+1}=\cdots=m_{e_{r-1}}=1$ for some $k.$

Since all the $e_j$'s are greater than $1,$ any $Q_i$ for $k+1 \le i \le e_{r-1}$ does not contain any $S$-vertices. Therefore, each $s_j$ for any $j \in [r-2]$ is in one of $Q_1, \dots, Q_k.$ Let $B_i$ be the set of $j$'s where $s_j$ is $Q_i,$ for any $i: 1 \le i \le k.$ We check that $Q_i$ is the union of $P_j$ for all $j \in B_i$ and $G$ is the union of $P_{r-1}$ and $\cup_{i=1}^k Q_i.$

For each $i: 1 \le i \le e_{r-1},$ let $\gamma_i$ be the cycle obtained by reading the $[d]$-vertex set of $Q_i$ on $\C_\tau$ in clockwise order. Because the $Q_i$'s have property (i), by Lemma \ref{char2prd} with $\eta = \sigma_{r-1}^{-1}$ and $\mu = \tau,$ we have that $\prod_{i=1}^{e_{r-1}} \gamma_i$ is the cycle decomposition of $\tau \sigma_{r-1}^{-1}$. Moreover, since $\gamma_{k+1},\dots,\gamma_{e_{r-1}}$ are cycles of length $1,$ we have $\displaystyle \prod_{i=1}^{k} \gamma_i = \tau \sigma_{r-1}^{-1}.$

Let $i$ be an integer with $1 \le i \le k.$ One sees that $Q_i \in \cG_S(m_i, \#B_i+1; (e_j)_{j \in B_i}).$ It is clear that $Q_i$ is a tree because $G$ is a tree. We then claim $Q_i$ also satisfies the following:
\begin{ilist}
\item[(ii)] Any $[d]$-vertex of $Q_i$ has CICPP on $(Q_i,\gamma_i).$
\end{ilist}
We can prove (ii) similarly as we did in the first half of this proof by using Lemma \ref{lem:iffCPP}. We omit the details.

Now by the induction hypothesis, we have that $Q_i \in \cG_S^*(m_i, \#B_i+1, \gamma_i; (e_j)_{j \in B_i}),$ which implies that $\prod_{j \in B_i} \sigma_j = \gamma_i.$
Since for any $j_1 \in B_{i_1}$ and $j_2 \in B_{i_2}$ with $i_1 \neq i_2$, we have that $\supp(\sigma_{j_1})$ and $\supp(\sigma_{j_2})$ are disjoint, $\sigma_{j_1}$ and $\sigma_{j_2}$ commute. Hence,
$$\prod_{j=1}^{r-2} \sigma_j = \prod_{i=1}^k \prod_{j \in B_i} \sigma_j = \prod_{i=1}^k \gamma_i = \tau \sigma_{r-1}^{-1}.$$
Therefore, $\sigma_1 \dots \sigma_{r-1} = \tau.$

Thus, we proved that the proposition holds for $r = r_0.$
\end{proof}

\section{Proof of Theorem \ref{thm:bijg2t}}\label{sec:proof}
Let $S = \{s_1 < s_2 < \cdots < s_{r-1}\}$ be a set of positive integers disjoint from $\{1,\dots,d\}.$ Also, by convention, we set $s_0 = 0.$ (So $s_0 < s_1 < \cdots < s_{r-1}.$)

We define $\LMR_S^*(1, e_1-1, \dots, e_{r-1}-1)$ to be the set of all labeled multi-noded rooted trees $\Phi^\L(G)$ associated to factorization graphs $G \in \cG_S^*(d,r,(1~2~\cdots ~d);$ $e_1,$ $\dots,$ $e_{r-1}).$ Then
$$\LMR_S^*(1, e_1-1, \dots, e_{r-1}-1) \subset \LMR_S(1, e_1-1, \dots, e_{r-1}-1).$$
The map $\Phi$ can be factored into two steps, as shown in Figure \ref{mapPhi}.

\begin{figure}[h,t]
\begin{center}
\begin{pspicture}(10pt,0pt)(260pt,50pt)
\psset{unit=10pt}
\rput[rt](10,5){$\cG^*_S(d,r,(1~2~\cdots ~d); e_1, \dots, e_{r-1})$}
\rput[lt](15,5){$\LMR_S^*(1, e_1-1, \dots, e_{r-1}-1)$}
\rput[lb](15,0){$\MR_S(1, e_1-1, \dots, e_{r-1}-1)$}

\psline{->}(11,4.5)(14.5,4.5)
\psline{->}(20,4)(20,1)
\psline{->}(5,4)(14.5,0.7)
\rput[bc](12.5,5){$\Phi^\L$}
\rput[bl](20.3,2){removing labels of nodes}
\rput[br](9,1.5){$\Phi$}
\end{pspicture}
\end{center}
\caption{}
\label{mapPhi}
\end{figure}

Hence, Theorem \ref{thm:bijg2t} follows from the following two lemmas.

\begin{lem}\label{lem:bijg2lt}
Suppose $\sum_{j=1}^{r-1} (e_j-1) = d-1$.
Then $\Phi^\L$ is a bijection from $\cG_S^*(d,r,(1~2~\cdots ~d);e_1,\dots,e_{r-1})$ to $\LMR_S^*(1,$ $e_1-1,$ $\dots,$ $e_{r-1}-1).$
\end{lem}

\begin{lem}\label{lem:bijlt2t}
Suppose $\sum_{j=1}^{r-1} (e_j-1) = d-1$.
For any $M \in \MR_S(1, e_1-1, \dots, e_{r-1}-1),$ there exists a unique labeling $\l$ of the nodes of $M$ such that $(M,\l) \in \LMR_S^*(1, e_1-1, \dots, e_{r-1}-1).$
\end{lem}

\begin{proof}[Proof of Lemma \ref{lem:bijg2lt}]
Given any $(M,\l) \in \LMR_S(1, e_1-1, \dots, e_{r-1}-1),$ we define $\Psi(M,\l)$ to be the $S$-$[d]$ bipartite graph $G$ whose edge set consists of $\{s, \nu\}$ for which $\nu$ is either a node contained in vertex $s$ in $M$ or the parent of $s$ in $M.$ It is clear that $G$ is in $\cG_S(d,r;e_1,\dots,e_{r-1})$ and is connected. Then by Lemma \ref{lem:tree}, $G$ is a tree. Hence, $\Psi(M,\l)$ is a map from $\LMR_S(1, e_1-1, \dots, e_{r-1}-1)$ to the set of $S$-$[d]$ bipartite trees.

For any $G \in \cG_S^*(d,r,(1,\dots,d);e_1,\dots,e_{r-1}),$ we have that $\Psi(\Phi^\L(G)) = G.$ Hence, $\Phi^\L$ is injective. Thus, the lemma follows.
\end{proof}

In order to prove Lemma \ref{lem:bijlt2t}, we need to discuss properties of the labeling $\l$ of any $(M, \l) \in \LMR_S^*(1, e_1-1, \dots, e_{r-1}-1).$ For convenience, we give the following definitions:
\begin{defn}
Given $(M, \l) \in \LMR_S(1, e_1-1, \dots, e_{r-1}-1),$ and any subgraph $M'$ of $M,$ we denote by $\l(M')$ the set of labels of the nodes in $M'.$

For any node $\nu$, we denote by $M_\nu$ the subtree of $M$ whose root has the single node $\nu.$

For any vertex $s$, we denote by $M_{s}$ the subtree of $M$ rooted at $s.$
\end{defn}

\begin{lem}\label{lem:charLMR}
Assume $\sum_{j=1}^{r-1} (e_j-1) = d-1$. Let $(M,\l) \in \LMR_S(1, e_1-1,$ $\dots,$ $e_{r-1}-1).$ Then $(M,\l) \in \LMR_S^*(1, e_1-1, \dots, e_{r-1}-1)$ if and only if there exist $1 \le \alpha_\nu \le \beta_\nu \le d$ for each node $\nu$ of $M$ and $1 \le \alpha_j' \le \beta_j' \le d$ for each vertex $s_j$ of $M$ satisfying the following conditions:
\begin{ilist}
\itm For any $\nu$ a node of $M,$ $\l(M_\nu) = [\alpha_\nu, \beta_\nu] := \{\alpha_\nu, \alpha_\nu+1, \dots, \beta_\nu\}.$

\itm For any $s_j$ a vertex of $M,$ $\l(M_{s_j}) = [\alpha_j', \beta_j'].$

\itm Suppose $\nu$ is a node contained in the vertex $s_j,$ and $s_{j_1}, \ldots, s_{j_\ell}$ are the vertices connected to $\nu$ with $j_1 < \cdots < j_k < j < j_{k+1} < \ldots < j_\ell$ for some $0 \leq k \leq \ell$. Then $\{\l(\nu)\}, [\alpha_{j_1}', \beta_{j_1}'], \dots, [\alpha_{j_\ell}', \beta_{j_\ell}']$ partition $[\alpha_\nu,\beta_\nu]$ into consecutive pieces with $\beta_{j_{k}}' <\cdots <  \beta_{j_{1}}' <  \l(\nu)  < \beta_{j_{\ell}}' < \cdots < \beta_{j_{k+1}}'.$

\itm Suppose $s_j$ is a vertex of $M.$ Let $\nu_1, \nu_2, \ldots, \nu_{e_j-1}$ be the nodes in $s_j$ from left to right. Then $[\alpha_{\nu_1}, \beta_{\nu_1}], \dots, [\alpha_{\nu_{e_j-1}}, \beta_{\nu_{e_j-1}}]$ partition $[\alpha_j',\beta_j']$ into consecutive pieces with $\beta_{\nu_1} < \cdots <  \beta_{\nu_{e_j-1}}.$

\end{ilist}
\end{lem}

\begin{ex}
Let $M$ be the labeled multi-noded rooted tree in Figure \ref{fig_MutiLabelled}.
We verify Lemma \ref{lem:charLMR} for some parts of $M.$ For vertex $s_9$ and nodes contained in $s_9,$ we have
\[ \l(M_{s_9}) =  [3, 12], \ \  \l(M_5) =  [3,5], \ \ \l(M_6) =  [6,6], \ \ \l(M_{11}) =  [7,11], \ \ \l(M_{12}) =  [12,12]. \]
Clearly, we have $\l(M_5), \l(M_6), \l(M_{11}), \l(M_{12})$ partition $\l(M_{s_9})$ into consecutive pieces with $\beta_5 = 5 < \beta_6=6 < \beta_{11} = 11 < \beta_{12} = 12.$ (Note that it is a coincidence that all $\beta_{\nu} = \nu$ for the four nodes we discussed. It is not always the case, e.g., $\beta_{15} = 18 \neq 15.$)

$s_1$ and $s_7$ are the vertices connected to the node $11.$ We have $1 < 7 < 9,$ and $\{11\},$ $\l(M_{s_1})=[10, 10]=\{10\}$, $\l(M_{s_7}) =[7,9]=\{7,8,9\}$ partition $\l(M_{11})=[7,11]$ into consecutive pieces. Furthermore, we have $\beta_{7}' = 9 < \beta_{1}' = 10 < 11.$
\end{ex}

\begin{proof}[Proof of Lemma \ref{lem:charLMR}]
Suppose $(M,\l) \in \LMR_S^*(1, e_1-1, \dots, e_{r-1}-1).$ Then $(M,\l) = \Phi^\L(G)$ for some $G \in \cG_S^*(d,r,(1~2~\cdots ~d); e_1,\dots, e_{r-1}).$ By Proposition \ref{prop:bijgraph}, $G$ satisfies (2) and (3) of Proposition \ref{prop:bijgraph}. It follows directly that $\l(M_\nu)$ for any node $\nu$ and $\l(M_{s_j})$ for any vertex $s_j$ are consecutive pieces on the circle $\C.$ Furthermore, since $1$ is the label of the node in the root, one sees each consecutive piece is actually a consecutive piece of $[1,d].$ Hence, we can define $\alpha_\nu, \beta_\nu$ and $\alpha_j', \beta_j'$ such that (i) and (ii) are satisfied.

Let $\nu$ be a node of $M.$ Assume $\nu$ is the single node labeled by $1$ in the root $s_0=0$. Because $s_0 < s_{j_1} < \cdots < s_{j_\ell},$ (iii) follows from the fact that the $[d]$-vertex $1$ has CICPP on $(G,(1~2\cdots~d))$. Suppose $\nu$ is not in the root. We denote by $\bar{M}_\nu$ the tree obtained from $M$ be removing $M_\nu.$
Then the fact that $\nu$ has CICCP on $(G,(1~2\cdots~d))$ implies that $\l(M_{s_{j_1}}), \dots, \l(M_{s_{j_k}}), \l(\bar{M}_\nu), \l(M_{s_{j_{k+1}}}), \dots, \l(M_{s_{j_\ell}})$ are consecutive pieces on $\C$ starting from $\nu$ in counterclockwise order. Note that $\l(\bar{M}_\nu)$ contains node $1,$ and the union of $\l(M_{s_{j_1}}), \dots, \l(M_{s_{j_k}}), \l(M_{s_{j_{k+1}}}),$ $\dots,$ $\l(M_{s_{j_\ell}})$ and $\{\l(\nu)\}$ is $\l(M_\nu).$ Thus, (iii) follows.

Let $j \in \{0,1,\dots, r-1\}.$ If $j=0,$ (iv) clearly holds. Suppose $j \in [r-1].$ One sees that $s_j$ having CPP on $(G,(1~2\cdots~d))$ implies that $\l(M_{\nu_1}), \dots, \l(M_{\nu_{e_{j}-1}})$ partition $\l(M_{s_j})$ into consecutive pieces. Furthermore, when we construct $(M,\l) = \Phi^\L(G)$ from $G,$ we require the labels of the nodes in $s_j$ to be in increasing order from left to right. It follows that $\beta_{\nu_1} < \cdots <  \beta_{\nu_{e_j-1}}.$ Therefore, (iv) holds.

Now we prove the other direction. Suppose $(M,\l) \in \LMR_S(1, e_1-1, \dots, e_{r-1}-1)$ and there exist $1 \le \alpha_\nu \le \beta_\nu \le d$ for each node $\nu$ of $M$ and $1 \le \alpha_j' \le \beta_j' \le d$ for each vertex $s_j$ of $M$ satisfying (i)-(iv).
Let $\Psi$ be the map from $\LMR_S(1, e_1-1, \dots, e_{r-1}-1)$ to the set of $S$-$[d]$ bipartite trees defined in the proof of Lemma \ref{lem:bijg2lt}, and define $G := \Psi(M,\l).$
We can reverse the proof in the last two paragraphs to show that (iii) and (iv) imply that $G$ satisfies (2) and (3) of Proposition \ref{prop:bijgraph}. Since $G$ is also a tree, using Proposition \ref{prop:bijgraph}, we conclude that $G \in \cG_S^*(d,r,(1~2~\cdots~d);e_1,\dots,e_{r-1}).$ It is sufficient to show that $\Phi^\L(G) = (M, \l).$ However, one checks that for any $(M,\l) \in \LMR_S(1, e_1-1, \dots, e_{r-1}-1)$, $\Phi^\L(\Psi(M, \l)) = (M, \l)$ if and only if the following two conditions hold:
\begin{enumerate}
\item The label of the single node in the root $s_0$ of $M$ is $1$.
\item For any $j\in [r-1]$ the labels of the nodes in $s_j$ are in increasing order from left to right.
\end{enumerate}
However, (1) follows from (iii) by letting $\nu$ be the single node in $s_0,$ and (2) follows from the condition $\beta_{\nu_1} < \cdots <  \beta_{\nu_{e_j-1}}$ in (iv).
\end{proof}

\begin{proof}[Proof of Lemma \ref{lem:bijlt2t}]
Let $M \in \MR_S(1, e_1-1, \dots, e_{r-1}-1).$
By Lemma \ref{lem:charLMR}, it is equivalent to prove that there exists a unique choice of a labeling $\l$ for the nodes of $M$ with set $[d],$ and integers $1 \le \alpha_\nu \le \beta_\nu \le d$ for each node $\nu$ of $M$ and integers $1 \le \alpha_j' \le \beta_j' \le d$ for each vertex $s_j$ of $M$ such that (i)-(iv) of Lemma \ref{lem:charLMR} are satisfied.

For any vertex $s_j$, we say it is a {\it level-$m$} vertex if it has distance $m$ to the root $s_0.$
We call a node a {\it level-$m$} node if it is inside a level-$m$ vertex.
We will describe an algorithm to choose the unique $\l$, $\alpha_\nu, \beta_\nu$ and $\alpha_j',\beta_j'.$
The algorithm will assign values in the order of levels: At step (0), we define $\alpha_0'$ and $\beta_0'$ for the root $s_0;$
at step (2m+1) (for $m \ge 0$), we define $\alpha_\nu$ and $\beta_\nu$ for all level-$m$ nodes;
at step (2m+2) (for $m \ge 0$), we define $\l(\nu)$ for all level-$m$ nodes, and define $\alpha_j'$ and $\beta_j'$ for all level-$(m+1)$ vertices.

\begin{enumerate}
\item[(0)] For the root $s_0$ of $M,$ since $M_{s_0} = M,$ the set of labels in $M_{s_0}$ is just $[d].$ Therefore, there is a unique way to choose $\alpha_0' = 1$ and $\beta_0' = d.$
\item[(2m+1)] Suppose for any vertex $s_j$ at level-$m,$ $\alpha_j'$ and $\beta_j'$ are defined.

Let $s_j$ be a vertex at level-$m$ and let $\nu_1, \nu_2, \ldots, \nu_{e_j-1}$ be the nodes in $s_j$ from left to right. Let $n_i$ be the number of nodes in $M_{\nu_i}$ for each $1 \le i \le e_j-1.$ Since $\alpha_j'$ and $\beta_j'$ are defined already, one sees that there is a unique way to choose $\alpha_{\nu_1}, \beta_{\nu_1}, \dots, \alpha_{\nu_{e_j-1}}, \beta_{\nu_{e_j-1}}$ such that (iv) of Lemma \ref{lem:charLMR} is satisfied for $s_j:$
$$\alpha_{\nu_i} := \alpha_j'+\sum_{t=1}^{i-1}n_t, \ \ \beta_{\nu_i} := \alpha_j' -1 + \sum_{t=1}^{i}n_t, \ \ \forall 1 \le i \le e_j-1.$$

Therefore, in this step, we can define $\alpha_\nu$ and $\beta_\nu$ for all the nodes at level-$m.$

\item[(2m+2)] Suppose for any vertex $\nu$ at level-$m,$  $\alpha_\nu$ and $\beta_\nu$ are defined.

 Let $\nu$ be a level-$m$ node contained in vertex $s_j,$ and $s_{j_1}, \ldots, s_{j_\ell}$ the vertices connected to $\nu$ with $j_1 < \cdots < j_k < j < j_{k+1} < \ldots < j_\ell$ for some $1 \leq k \leq \ell$. Clearly, $s_{j_1},\dots, s_{j_\ell}$ are level-$(m+1)$ vertices. Let $n_1,n_2,\dots, n_\ell$ be the number of nodes in $M_{j_1}, M_{j_2}, \dots, M_{j_\ell}.$ Since $\alpha_\nu$ and $\beta_\nu$ are defined already, one sees that there is a unique way to choose $\l(\nu),\alpha_{j_1}', \beta_{j_1}', \dots, \alpha_{j_{e_j-1}}', \beta_{j_{e_j-1}}'$ such that (iii) of Lemma \ref{lem:charLMR} is satisfied for $\nu:$
 $$\l(\nu) = \alpha_\nu + \sum_{t=1}^{k}n_t;$$
 $$\alpha_{j_i}' := \alpha_\nu+\sum_{t=i+1}^{k} n_t, \ \ \beta_{j_i}' := \alpha_\nu -1 + \sum_{t=i}^{k}n_t, \ \ \forall 1 \le i \le k;$$
 $$\alpha_{j_i}' := \alpha_\nu+\sum_{t=1}^{k}n_t+1+\sum_{t=i+1}^{\ell} n_t, \ \beta_{j_i}' := \alpha_\nu + \sum_{t=1}^{k}n_t + \sum_{t=i}^{\ell}n_t, \ \forall k+1 \le i \le \ell.$$

Therefore, in this step, we define labels for all the nodes at level-$m$ and $\alpha_j'$ and $\beta_j'$ for all the vertices at level-$(m+1).$
\end{enumerate}
It is easy to see that this algorithm defines the unique solution to $\l, \alpha_\nu, \beta_\nu, \alpha_j', \beta_j'$ that satisfies (i)-(iv) of Lemma \ref{lem:bijlt2t}.
\end{proof}

We proved Lemma \ref{lem:bijg2lt} and \ref{lem:bijlt2t}. Hence, Theorem \ref{thm:bijg2t} follows.


\begin{thebibliography}{99}
\small \setlength{\itemsep}{-.8mm}

\bibitem{Denes}
J. D\'{e}nes, {\em The representation of a permutation as the product of
a minimal number of transpositions and its connection with the theory of
graphs}, Publ. Math. Institute Hung. Acad. Sci. {\bf 4} (1959), 63--70.

\bibitem{GouldenJacksonEJC} I.P. Goulden and D.M. Jackson, {\em The combinatorial relationship between trees, cacti and certain connection coefficients for the symmetric group}, European J. Combin. {\bf 13} (1992), 357--365.

\bibitem{GouldenJakson} I.P. Goulden and D.M. Jackson, {\em Transitive factorizations into transpositions and holomorphic mappings on the sphere}, Proc. Amer.
Math. Soc. {\bf 125} (1997), 51--60.

\bibitem{GouldenPepper} I.P. Goulden and S. Pepper, {\em Labelled trees and factorizations of a cycle into transpositions}, Discrete Math. {\bf 113} (1993), 263--268.

\bibitem{Treelike} I.P. Goulden and A. Yong, {\em Tree-like properties of cycle factorizations}, J. Combin. Theory Ser. A {\bf 98}(1) (2002), 106--117.

\bibitem{hurwitz} A. Hurwitz, {\em Ueber Riemann'sche Fl\"{a}chen mit gegebenen Verzweigungspunkten}, Mathematische Annalen {\bf 39} (1891), 1--60.

\bibitem{Irving} J. Irving, {\em Minimal transitive factorizations of permutations into cycles}, Canad. J. Math. {\bf 61} (2009), 1092--1117.

\bibitem{LandoZvonkine} S. Lando and D. Zvonkine, {\em On the multiplicities of the Lyashko-Looijenga map on the strata of the discriminant}, Funkt. Anal. Appl. {\bf 33} (3) (1999), 21--34.

\bibitem{LiuOsserman} F. Liu and B. Osserman, {\em The irreducibility of certain pure-cycle Hurwitz spaces}, Amer. J. Math. {\bf 130} (2008), 1687--1708.

\bibitem{Moszkowski} P. Moszkowski, {\em A solution to a problem of D\'{e}nes: a bijection
between trees and factorizations of cyclic permutations}, European
J. Combinatorics {\bf 10} (1989), 13--16.

\bibitem{Springer}
C. M. Springer, {\em Factorizations, trees, and cacti}, Eighth International Conference on Formal Power Series and Algebraic Combinatorics, University of Minnesota, June 25-29, 1996, 427--438.

\bibitem{stanleyec1} R.P. Stanley, {\em Enumerative Combinatorics}, vol.~1, Cambridge Studies in Advanced Mathematics, vol. 49, Cambridge University Press, Cambridge, 1997.

\bibitem{stanleyec2} R.P. Stanley, {\em Enumerative Combinatorics}, vol.~2, Cambridge Studies in Advanced Mathematics, vol. 62, Cambridge University Press, Cambridge, 1999.

\end{thebibliography}
\end{document}